%&amstex          
\input amstex\documentstyle{amsppt}  
\pagewidth{12.5cm}\pageheight{19cm}\magnification\magstep1
\topmatter
\title On conjugacy classes in a reductive group\endtitle
\author G. Lusztig\endauthor
\address{Department of Mathematics, M.I.T., Cambridge, MA 02139}\endaddress
\thanks{Supported in part by National Science Foundation grant DMS-1303060.}\endthanks
\endtopmatter   
\document
\define\bx{\boxed}
\define\rk{\text{\rm rk}}

\define\Irr{\text{\rm Irr}}

\define\tcs{\ti{\cs}}

\define\ucl{\un{cl}}

\define\ul{\un l}

\define\mpb{\medpagebreak}

\define\bR{\bar R}

\define\si{\sim}

\define\sqc{\sqcup}

\define\lb{\linebreak}

\define\op{\oplus}
   
\define\part{\partial}
\define\emp{\emptyset}
\define\imp{\implies}
\define\ra{\rangle}
\define\n{\notin}

\define\m{\mapsto}
\define\do{\dots}
\define\la{\langle}

\define\lra{\leftrightarrow}

\define\sub{\subset}    
\define\bxt{\boxtimes}
\define\T{\times}
\define\ti{\tilde}
\define\nl{\newline}
\redefine\i{^{-1}}

\define\un{\underline}

\define\ot{\otimes}
\define\bbq{\bar{\QQ}_l}

\define\Hom{\text{\rm Hom}}

\define\Aut{\text{\rm Aut}}
\define\Ind{\text{\rm Ind}}

\define\tr{\text{\rm tr}}

\define\a{\alpha}
\redefine\b{\beta}

\define\g{\gamma}
\redefine\d{\delta}
\define\e{\epsilon}

\define\io{\iota}

\define\ph{\phi}
\define\ps{\psi}
\define\r{\rho}
\define\s{\sigma}
\redefine\t{\tau}
\define\th{\theta}

\redefine\l{\lambda}

\define\x{\xi}

\redefine\G{\Gamma}

\define\Ph{\Phi}

\define\kk{\bold k}

\define\rr{\bold r}

\define\CC{\bold C}

\define\FF{\bold F}

\define\HH{\bold H}

\define\NN{\bold N}

\define\QQ{\bold Q}

\define\ZZ{\bold Z}

\define\cb{\Cal B}
\define\cc{\Cal C}

\define\co{\Cal O}
\define\cp{\Cal P}

\define\cs{\Cal S}
\define\ct{\Cal T}
\define\cu{\Cal U}

\define\fb{\frak b}

\define\fg{\frak g}

\define\fp{\frak p}

\define\fz{\frak z}

\define\fU{\frak U}

\define\tz{\ti z}

\define\tA{\ti A}

\define\tE{\ti E}

\define\sha{\sharp}

\define\che{\check}

\define\BO{Bo}
\define\CA{C}
\define\CAR{Ca}
\define\DL{DL}
\define\GH{CH}
\define\GP{GP}
\define\KL{KL}
\define\LCL{L1}
\define\ORA{L2}
\define\ICC{L3}
\define\CDG{L4}
\define\UNI{L5}
\define\SPEC{L6}
\define\REM{L7}
\define\WEU{L8}
\define\CSM{L9}
\define\WEUII{L10}
\define\WEUIII{L11}
\define\DIST{L12}
\define\UNAC{L13}
\define\LS{LS1}
\define\LSS{LS2}
\define\PE{Pe}
\define\SPA{Spa}
\define\SPR{Spr}
\define\ST{St}

\head Introduction\endhead
\subhead 0.1\endsubhead
Let $\kk$ be an algebraically closed field of characteristic $p\ge0$ and let $G$ be a connected reductive 
algebraic group over $\kk$. Let $W$ be the Weyl group of $G$. Let $cl(W)$ be the set of conjugacy classes of
$W$. 

In \cite{\ST} Steinberg has defined the notion of regular element in $G$ (an element whose conjugacy class
has dimension as large as possible, that is $\dim(G)-\rk(G)$) and showed that the set of regular elements
in $G$ form an open dense subset $G_{reg}$. The goal of this paper is to define a partition of $G$ into 
finitely many strata, one of which is $G_{reg}$. Each stratum of $G$ is a union of 
conjugacy classes of $G$ of the same dimension. The set of strata is naturally indexed by a set
which depends only on $W$ as a Coxeter group, not on the underlying root system and not on the ground field
$\kk$. We give two descriptions of the indexing set above: 

(i) one in terms of a class of irreducible representations of $W$ which we call $2$-special representations
(they are obtained by truncated induction from special representations of certain reflection subgroups of 
$W$); 

(ii) one in terms of $cl(W)$ (modulo a certain equivalence relation).
\nl
In the case where $W$ is irreducible we give a third description of the indexing set above:

(iii) in terms of the sets of unipotent classes in the various versions of $G$ over $\bar{\FF}_r$ for a 
variable prime number $r$, glued together according to the set of unipotent classes in the version of $G$ 
over $\CC$. 
\nl
The definition of strata in the form (i) and (iii) are based on Springer's correspondence (see \cite{\SPR} 
when $p=0$ or $p\gg0$ and \cite{\ICC} for any $p$) connecting irreducible representations of $W$ with 
unipotent classes; when $W$ is irreducible, the definition of strata in the form (iii) is related to that 
in the form (ii) by the results of \cite{\WEU,\WEUII} connecting $cl(W)$ with unipotent classes in $G$. 

Since (i),(ii) are two incarnations of our indexing set, they are in canonical
bijection with each other. In particular we obtain a canonical map from $cl(W)$ to the set of irreducible
representations of $W$ whose image consists of the $2$-special representations (when $G$ is $GL_n(\kk)$ this
is a bijection).
We also show that the dimension of a conjugacy class in a stratum  of $G$ is independent of the ground
field. (This statement makes sense since the parametrization of the strata is independent of the ground 
field.) In particular, we see that if $n\ge1$, then the following three conditions on an integer $k$ are 
equivalent:

-there exists a conjugacy class of dimension $k$ in $SO_{2n+1}(\CC)$;

-there exists a conjugacy class of dimension $k$ in $Sp_{2n}(\CC)$;

-there exists a conjugacy class of dimension $k$ in $Sp_{2n}(\bar F_2)$.
\nl
The proof shows that the following fourth condition is equivalent to the three conditions above: there 
exists a unipotent conjugacy class of dimension $k$ in $Sp_{2n}(\bar F_2)$.

In \S5 we sketch an alternative approach to the definition of strata which is based on an extension
of the ideas in \cite{\WEU}, and Springer's correspondence does not appear in it.

In \S6 we dicuss extensions of our results to the Lie algebra of $G$ and to the case where $G$ is replaced
by a disconnected reductive group. We also define a partition of the set of compact regular semisimple 
elements in a loop group into strata analogous to the partition of $G$ into strata. Moreover we give a 
conjectural description of the strata of $G$ (assuming that $\kk=\CC$) which is based on an extension of a
construction in \cite{\KL}.

\subhead 0.2\endsubhead 
{\it Notation.}
For an algebraic group $H$ over $\kk$ we denote by $H^0$ the identity component of $H$. For a subgroup $T$ 
of $H$ we denote by $N_HT$ the normalizer of $T$ in $H$. Let $\fg$ be the Lie algebra of $G$.
For $g\in G$ we denote by $Z_G(g)$ the centralizer of $g$ in $G$ and by $g_s$ (resp. $g_u$) the semisimple 
(resp. unipotent) part of $g$. Let $\cb$ be the variety of Borel subgroups of $G$.
Let $\cb_g=\{B\in\cb;g\in B\}$. Let $l$ be a prime number $\ne p$. For an 
algebraic variety $X$ over $\kk$ we denote by $H^i(X)$ the $l$-adic cohomology of $X$ in degree $i$; if $X$ 
is projective let $H_i(X)=\Hom(H^i(X),\QQ_l)$.

For any (finite) Weyl group $\G$ we denote by $\Irr\G$ a set of representatives for the isomorphism classes
of irreducible representations of $\G$ over $\QQ$. For any $\t\in\Irr W$ let $n_\t$ be the smallest integer 
$i\ge0$ such that $\t$ appears with $>0$ multiplicity in the $i$-th symmetric power of the reflection 
representation of $W$; if this multiplicity is $1$ we say that $\t$ is {\it good}.

A {\it bipartition} is a sequence $\l=(\l_1,\l_2,\l_3,\do)$ in $\NN$ such that $\l_m=0$ for $m\gg0$ and 
$\l_1\ge\l_3\ge\l_5\ge\do$, $\l_2\ge\l_4\ge\l_6\ge\do$. We write $|\l|=\l_1+\l_2+\l_3+\do$. We say that $\l$
is a bipartition of $n$ if $|\l|=n$. Let $BP^n$ be the set of bipartitions of $n$. Let $e,e'\in\NN$. We say 
that a bipartition $(\l_1,\l_2,\l_3,\do)$ has excess $(e,e')$ if $\l_i+e\ge\l_{i+1}$ for $i=1,3,5,\do$ and 
$\l_i+e'\ge\l_{i+1}$ for $i=2,4,6,\do$. Let $BP^n_{e,e'}$ be the set of bipartitions of $n$ which have 
excess $(e,e')$. A {\it partition} is a sequence $\l=(\l_1,\l_2,\l_3,\do)$ in $\NN$ such that $\l_m=0$ for 
$m\gg0$ and $\l_1\ge\l_2\ge\l_3\ge\do$. Thus a partition is the same as a bipartition of excess $(0,0)$. On 
the other hand, a bipartition is the same as an ordered pair of partitions 
$((\l_1,\l_3,\l_5,\do),(\l_2,\l_4,\l_6,\do))$.

Let $\cp=\{2,3,5,\do\}$ be the set of prime numbers.

\head 1. The $2$-special representations of a Weyl group\endhead
\subhead 1.1\endsubhead
Let $V,V^*$ be finite dimensional $\QQ$-vector spaces
with a given perfect bilinear pairing $\la,\ra:V\T V^*@>>>\QQ$.
Let $R$ (resp. $\che R$) be a finite subset of $V-\{0\}$ (resp. $V^*-\{0\}$) with a given bijection
$\a\lra\che\a$, $R\lra\che R$, such that $\la\a,\che\a\ra=2$ for any $\a\in R$ and 
$\la\a,\che\b\ra\in\ZZ$ for any $\a,\b\in R$; it is assumed that 
$\b-\la\b,\che\a\ra\a\in R$, $\che\b-\la\a,\che\b\ra\che\a\in\che R$ for any $\a,\b\in R$ and that
$\a\in R\imp\a/2\n R$. Thus, $(V,V^*,R,\che R)$ is a reduced root system.
Let $V_0$ (resp. $V_0^*$) be the $\QQ$-subspace of $V$ (resp. $V^*$) spanned by $R$ (resp.
$\che R$). Let $\rk(R)=\dim V_0=\dim V_0^*$. Let $W$ be the
(finite) subgroup of $GL(V)$ generated by the reflections $s_\a:x\m x-\la x,\che\a\ra\a$ in $V$ for
various $a\in R$; it may be identified with the subgroup of $GL(V^*)$ generated by the reflections 
${}^ts_a:x'\m x'-\la\a,x'\ra\che\a$ in $V^*$ for various $\a\in R$. 
For any $e\in V$ let $R_e=\{\a\in R;\la e,\che\a\ra\in\ZZ\}$, $\che R_e=\{\che\a;\a\in R_e\}$; note that 
$(V,V^*,R_e,\che R_e)$ is a root system with Weyl group $W_e=\{w\in W;w(e)-e\in\sum_{\a\in R}\ZZ\a\}$.
Similarly, for any $e'\in V^*$ let $R_{e'}=\{\a\in R;\la \a,e'\ra\in\ZZ\}$, 
$\che R_{e'}=\{\che\a;\a\in R_{e'}\}$; note that $(V,V^*,R_{e'},\che R_{e'})$ is a root system with
Weyl group $W_{e'}=\{w\in W;w(e')-e'\in\sum_{\a\in R}\ZZ\che\a\}$.
For any $(e,e')\in V\T V^*$ let $R_{e,e'}=R_e\cap R_{e'}$, $\che R_{e,e'}=\che R_e\cap\che R_{e'}$. Then
$(V,V^*,R_{e,e'},\che R_{e,e'})$ is a root system; let $W_{e,e'}$ be its Weyl group (a subgroup of
$W_e\cap W_{e'}$). Note that $W_{0,e'}=W_{e'}$, $W_{e,0}=W_e$, $W_{0,0}=W$.
For $E\in\Irr(W_{e,e'})$ let $n_E$ be as in 0.2. 

Let $(e_1,e'_1)\in V\T V^*$, $(e_2,e'_2)\in V\T V^*$ be such that
$R_{e_1,e'_1}\sub R_{e_2,e'_2}$ (so that $W_{e_1,e'_1}\sub W_{e_2,e'_2}$). In this case, if
$E\in\Irr(W_{e_1,e'_1})$ is good, there is a unique $E_0\in\Irr(W_{e_2,e'_2})$ such that $E_0$ appears in 
$\Ind_{W_{e_1,e'_1}}^{W_{e_2,e'_2}}(E)$ and $n_{E_0}=n_E$, see \cite{\LS, 3.2}; moreover, $E_0$ is good. We 
set $E_0=j_{W_{e_1,e'_1}}^{W_{e_2,e'_2}}(E)$. Note that if we have also $R_{e_2,e'_2}\sub R_{e_3,e'_3}$ where
$(e_3,e'_3)\in V\T V^*$ then we have the transitivity property:
$$j_{W_{e_1,e'_1}}^{W_{e_3,e'_3}}(E)=j_{W_{e_2,e'_2}}^{W_{e_3,e'_3}}(j_{W_{e_1,e'_1}}^{W_{e_2,e'_2}}(E)).
\tag a$$
Let $\cs(W_{e,e'})\sub\Irr(W_{e,e'})$ be the set of {\it special} representations of $W_{e,e'}$, see 
\cite{\LCL}; note that any $E\in\cs(W_{e,e'})$ is good. Hence $j_{W_{e,e'}}^W(E)\in\Irr(W)$ is defined. We 
say that $E_0\in\Irr(W)$ is {\it $2$-special} if $E_0=j_{W_{e,e'}}^W(E)$ for some $(e,e')\in V\T V^*$ and
some $E\in\cs(W_{e,e'})$. Let $\cs_2(W)$ be the set of all $2$-special representations of $W$ (up to 
isomorphism). From the definition we see that

(b) {\it $\cs_2(W)$ is unchanged when $(V,V^*,R,\che R)$ is replaced by $(V^*,V,\che R,R)$.}
\nl
Let $\cs_1(W)$ (resp. ${}'\cs_1(W)$) be the set of all $E_0\in\Irr(W)$ such that $E_0=j_{W_e}^W(E)$ (resp. 
$E_0=j_{W_{e'}}^W(E)$) for some $e\in V$, $E\in\cs(W_e)$ (resp. $e'\in V^*$, $E\in\cs(W_{e'}$). 
The analogue of (b) with $\cs_2(W)$ replaced by $\cs_1(W)$ is not true in general; instead, if
$(V,V^*,R,\che R)$ is replaced by $(V^*,V,\che R,R)$ then $\cs_1(W)$ becomes ${}'\cs_1(W)$ and ${}'\cs_1(W)$
becomes $\cs_1(W)$. 

Now, for any 
$e'\in V^*$ the subset $\cs_1(W_{e'})\sub\Irr(W_{e'})$ is defined; it consists of all $E'\in\Irr(W_{e'})$ 
such that $E'=j_{W_{e,e'}}^{W{e'}}(E)$ for some $e\in V$ and some $E\in\cs(W_{e,e'})$. Note that any 
$E'\in\cs_1(W_{e'})$ is good. From (a) we see that

(c) {\it $\cs_2(W)$ consists of all $E_0\in\Irr(W)$ such that $E_0=j_{W_{e'}}^W(E')$ for some $e'\in V^*$ and
some $E'\in\cs_1(W_{e'})$.}
\nl
We say that $e'\in V^*$ (resp. $(e,e')\in V\T V^*$) is {\it isolated} if $\rk(R_{e'})=\rk(R)$ (resp. 
$\rk(R_{e,e'})=\rk(R)$). We show:

(d) {\it $\cs_2(W)$ consists of all $E_0\in\Irr(W)$ such that $E_0=j_{W_{e,e'}}^W(E)$ for some isolated 
$(e,e')\in V\T V^*$ and some $E\in\cs(W_{e,e'})$.}
\nl
Let $E_0\in\cs_2(W)$. By definition, we can find $(e,e')\in V\T V^*$ and $E\in\cs(W_{e,e'})$ such that 
$E_0=j_{W_{e,e'}}^W(E)$. We can find an isolated $e'_1\in V^*$ such that $R_{e'}$ is rationally closed in 
$R_{e'_1}$ that is, $R_{e'_1}\cap\sum_{\a\in R_{e'}}\QQ\a=R_{e'}$. Applying the analogous statement to 
$(V^*,V,\che R_{e'_1},R_{e'_1})$, $e$, instead of $(V,V^*,R,\che R)$, $e'$, we can find $e_1\in V$ such that
$\rk(R_{e_1}\cap R_{e'_1})=\rk(R_{e'_1})$ and $R_e\cap R_{e'_1}$ is rationally closed in 
$R_{e_1}\cap R_{e'_1}$. It follows that $(e_1,e'_1)$ is isolated and $R_e\cap R_{e'}$ is rationally closed 
in $R_{e_1}\cap R_{e'_1}$ hence $E_1:=j_{W_{e,e'}}^{W_{e_1,e'_1}}(E)$ is in $\cs(W_{e_1,e'_1})$, see 
\cite{\LCL}. By (a), we have $E_0=j_{W_{e_1,e'_1}}^W(E_1)$. This proves (d).

We have the following variant of (d):

(e) {\it $\cs_2(W)$ consists of all $E_0\in\Irr(W)$ such that $E_0=j_{W_{e'}}^W(\tE)$ for some isolated
$e'\in V^*$ and some $\tE\in\cs_1(W_{e'})$.}
\nl
Let $E_0\in\cs_2(W)$. Let $E,e,e'$ be as in (d). We have $E=j_{W_{e'}}^W(\tE)$ where 
$\tE=j_{W_{e,e'}}^{W_{e'}}(E)\in\cs_1(W_{e'})$ and $\rk(R_{e'})=\rk(R)$. Conversely, if $e'\in V^*$ and 
$\tE\in\cs_1(W_{e'})$ then, by (c), $j_{W_{e'}}^W(\tE)\in\cs_2(W)$ (even without the assumption that 
$\rk(R_{e'})=\rk(R)$). This proves (e).

\mpb

Let $R'\sub R$ be such that (if $\che R'$ is the image of $R'$ under $R\lra\che R$),
$(V,V^*,R',\che R')$ is a root system (with Weyl group $W'$) and $R'$ is rationally closed in $R$. Note 
that $R'=R_e$ for some $e\in V$ and $R'=R_{e'}$ for some $e'\in V^*$. We show:

(f) If $E\in\cs_1(W')$ then $j_{W'}^W(E)\in\cs_1(W)$. 

(g) If $E\in\cs_2(W')$ then $j_{W'}^W(E)\in\cs_2(W)$. 
\nl
We prove (f). Let $e'\in V^*$ be such that $R'=R_{e'}$. We have $E=j_{W_{e,e'}}^{W_{e'}}(E')$ for some
$e\in V$ and some $E'\in\cs(W_{e,e'})$. Hence $j_{W'}^W(E)=j_{W_{e,e'}}^W(E')=j_{W_e}^W(E'')$ where 
$E''=j_{W_{e,e'}}^{W_e}(E')$. Now $R_{e,e'}$ is rationally closed in $R_e$ hence $E''\in\cs(W_e)$, see 
\cite{\LCL}. We see that $j_{W'}^W(E)\in\cs_1(W)$.

We prove (g). Let $e\in V$ be such that $R'=R_e$. We have $E=j_{W_{e,e'}}^{W_e}(E')$ for some $e'\in V^*$ 
and some $E'\in\cs_1(W_{e,e'})$. Hence $j_{W'}^W(E)=j_{W_{e,e'}}^W(E')=j_{W_{e'}}^W(E'')$ where 
$E''=j_{W_{e,e'}}^{W_{e'}}(E')$. Now $R_{e,e'}$ is rationally closed in $R_{e'}$ hence $E''\in\cs(W_{e'})$, 
see (f). We see that $j_{W'}^W(E)\in\cs_2(W)$.

\subhead 1.2\endsubhead
There  are unique direct sum decompositions $V_0=\op_{i\in I}V_i$, $V_0^*=\op_{i\in I}V_i^*$ such that
$R=\sqc_{i\in I}(R\cap V_i)$, $\che R=\sqc_{i\in I}(\che R\cap V_i)$ and for any $i\in I$,
$(V_i,V_i^*,R\cap V_i,\che R\cap V_i)$ is an irreducible root system for (with Weyl group $W_i$); the 
bijection $R\cap V_i\lra\che R\cap V_i$ is induced by $R\lra\che R$). We have canonically 
$W=\prod_{I\in I}W_i$ and $\cs_2(W)=\prod_{i\in I}\cs_2(W_i)$ (via external tensor product).

\subhead 1.3\endsubhead
In this subsection we assume that $(V,V^*,R,\che R)$ is irreducible. Now $W$ acts naturally on the set of 
subgroups $W'$ of $W$ of form $W_{e'}$ for various isolated $e'\in V^*$. The types of various $W'$ which 
appear in this way are well known and are described below in each case.

(a) $R$ of type $A_n$, $n\ge0$: $W'$ of type $A_n$.

(b) $R$ of type $B_n$, $n\ge2$: $W'$ of type $B_a\T D_b$ where $a\in\NN$, $b\in\NN-\{1\}$, $a+b=n$.

(c) $R$ of type $C_n$, $n\ge2$: $W'$ of type $C_a\T C_b$ where $a,b\in\NN$, $a+b=n$.

(d) $R$ of type $D_n$, $n\ge4$: $W'$ of type $D_a\T D_b$ where $a,b\in\NN-\{1\}$, $a+b=n$.

(e) $R$ of type $E_6$: $W'$ of type $E_6,A_5A_1,A_2A_2A_2$.

(f) $R$ of type $E_7$: $W'$ of type $E_7,D_6A_1,A_7,A_5A_2,A_3A_3A_1$.

(g) $R$ of type $E_8$: $W'$ of type $E_8,E_7A_1,E_6A_2,D_5A_3,A_4A_4$,$A_5A_2A_1,A_7A_1,A_8,D_8$.

(h) $R$ of type $F_4$: $W'$ of type $F_4,B_3A_1,A_2A_2,A_3A_1,B_4$.

(i) $R$ of type $G_2$: $W'$ of type $G_2,A_2,A_1A_1$.

(We use the convention that a Weyl group of type $B_n$ or $D_n$ with $n=0$ is $\{1\}$.)

\subhead 1.4\endsubhead
In this subsection we assume that $(V,V^*,R,\che R)$ is irreducible. Now $W$ acts naturally on the set of 
subgroups $W'$ of $W$ of form $W_{e,e'}$ for various isolated $(e,e')\in V\T V^*$. The types of various $W'$
which appear in this way are described below in each case. (For type $F_4$ and $G_2$ we denote by $\t$ a 
non-inner involution of $W$).

(a) $R$ of type $A_n$: $W'$ of type $A_n$.

(b) $R$ of type $B_n$ or $C_n$: $W'$ of type $B_a\T B_b\T D_c\T D_d$ where $a,b\in\NN$, $c,d\in\NN-\{1\}$,
$a+b+c+d=n$.

(c) $R$ of type $D_n$: $W'$ of type $D_a\T D_b\T D_c\T D_d$ where $a,b,c,d\in\NN-\{1\}$, $a+b+c+d=n$.

(d) $R$ of type $E_6$: $W'$ as in 1.3(e).

(e) $R$ of type $E_7$: $W'$ as in 1.3(f) and also $W'$ of type $D_4A_1A_1A_1$.

(f) $R$ of type $E_8$: $W'$ as in 1.3(g) and also $W'$ of type $D_6D_2,D_4D_4,A_3A_3A_1A_1$, $A_2A_2A_2A_2$.

(g) $R$ of type $F_4$: $W'$ as in 1.3(h), the images under $\t$ of the subgroups $W'$ of type $A_3A_1,B_4$ 
in 1.3(h) and also $W'$ of type $B_2B_2$.

(h) $R$ of type $G_2$: $W'$ as in 1.3(i) and the image under $\t$ of the subgroup $W'$ of type $A_2$ in 
1.3(i).

\subhead 1.5\endsubhead
If $R'\sub R$, $\che R'\sub\che R$ are such that $(V,V^*,R',\che R')$ is a root system (with the
bijection $R'\lra\che R'$ being induced by $R\lra\che R$) then, setting 
$\bR'=R\cap\sum_{\a\in R'}\QQ\a$, $\che{\bR'}=\che R\cap\sum_{\a\in R'}\QQ\che\a$, we obtain a root system 
$(V,V^*,\bR',\che{\bR'})$. We set 
$$N_{R'}=\sha(\sum_{\a\in\bR'}\ZZ\a/\sum_{\a\in R'}\ZZ\a)\in\ZZ_{\ge1}.$$
For any $e'\in V^*$ we set $N_{e'}=N_{R_{e'}}$.

Now let $r\in\cp$. Let $\cs^r_2(W)$ be the set of all $E_0\in\Irr(W)$ such that for some 
isolated $e'\in V^*$ with $N_{e'}=r^k$ for some $k\in\NN$ and for some $E\in\cs_1(W_{e'})$ we have 
$E_0=j_{W_{e'}}^W(E)$. Note that $\cs^1(W)\sub\cs^r_2(W)\sub\cs_2(W)$.

Now assume that $(V,V^*,R,\che R)$ is irreducible. We show:

(a) If $R$ is of type $A_n$, $n\ge0$ then $\cs^r_2(W)=\cs_2(W)=\cs_1(W)=\cs(W)$.

(b) If $R$ is of type $B_n$ or $C_n$, $n\ge2$ then $\cs^r_2(W)=\cs_1(W)$ if $r\ne2$ and 
$\cs^2_2(W)=\cs_2(W)$.

(c) If $R$ is of type $D_n$, $n\ge4$ then $\cs^r_2(W)=\cs_1(W)$ if $r\ne2$ and $\cs^2_2(W)=\cs_2(W)$.

(d) If $R$ is of type $E_6$ then $\cs^r_2(W)=\cs_2(W)=\cs_1(W)$.

(e) If $R$ is of type $E_7$ then $\cs^r_2(W)=\cs_1(W)$ if $r\ne2$ and $\cs^2_2(W)=\cs_2(W)$.

(f) If $R$ is of type $E_8$ then $\cs^r_2(W)=\cs_1(W)$ if $r\n\{2,3\}$ and 
$\cs^2_2(W)\cup\cs^3_2(W)=\cs_2(W)$.

(g) If $R$ is of type $F_4$ then $\cs^r_2(W)=\cs_1(W)$ if $r\ne2$ and $\cs^2_2(W)=\cs_2(W)$.

(h) If $R$ is of type $G_2$ then $\cs^r_2(W)=\cs_1(W)$ if $r\ne3$ and $\cs^3_2(W)=\cs_2(W)$.
\nl
We prove (a). In this case for any isolated $e'\in V^*$ we have $N_{e'}=1$ and the result follows from
1.1(d),(e), 1.3.

We prove (b),(c). In these cases for any isolated $e'\in V^*$, $N_{e'}$ is a power of $2$ (see 1.3) and the
equality $\cs^2_2(W)=\cs_2(W)$ follows from 1.1(e). Moreover, if $e'$ is isolated and $N_{e'}$ is not
divisible by $2$ then $W_{e'}=W$ so that for $r\ne2$ we have $\cs^r_2(W)=\cs_1(W)$.

In cases (d),(e),(f) we shall use the fact that for any $e'\in V^*$:

(i) we can find $e\in V$ such that $W_{e'}=W_e$, so that if $E\in\cs(W_{e'})$ then 
$j_{W_{e'}}^W(E)\in\cs_1(W)$.
\nl
(This property does not always hold in cases (g),(h).)

We prove (d). If $e'\in V^*$ is isolated and $W_{e'}\ne W$ then from 1.3 we see that $W_{e'}$ is of type 
$A_2A_2A_2$ or $A_5A_1$ so that $\cs_1(W_{e'})=\cs(W_{e'})$; using this and 1.1(e) we see that 
$\cs_2(W)=\cs^r_2=\cs_1(W)$. (We have used (i).)

We prove (e). If $e'\in V^*$ is isolated and $W_{e'}$ is not of type $E_7$ (with $N_{e'}=1$) or $D_6A_1$ 
(with $N_{e'}=2$) then from 1.3 we see that $W_{e'}$ is of type $A_7$ or $A_5A_2$ or $A_3A_3A_1$ so that 
$\cs_1(W_{e'})=\cs(W_{e'})$. We see that $\cs^r_2(W)=\cs_1(W)$ if $r\ne2$ and $\cs^2_2(W)=\cs_2(W)$. (We 
have used (i).)

We prove (f). If $e'\in V^*$ is isolated and $W_{e'}$ is not of type $E_8$ (with $N_{e'}=1$) or $E_7A_1$ 
(with $N_{e'}=2$) or $E_6A_2$ (with $N_{e'}=3$) or $D_5A_3$ (with $N_{e'}=4$) or $D_8$ (with $N_{e'}=2$) then
from 1.3 we see that $W_{e'}$ is of type $A_4A_4$ or $A_5A_2A_1$ or $A_7A_1$ or $A_8$, so that 
$\cs_1(W_{e'})=\cs(W_{e'})$; we see that $\cs^r_2(W)=\cs_1(W)$ if $r\n\{2,3\}$ and 
$\cs^2_2(W)\cup\cs^3_2(W)=\cs_2(W)$. (We have used (i).)

We prove (g). If $e'\in V^*$ is isolated and $W_{e'}$ is not of type $F_4$ (when $N_{e'}=1$) or $B_3A_1$ 
(with $N_{e'}$ a power of $2$) or $B_4$ (with $N_{e'}$ a power of $2$) then from 1.3 we see that $W_{e'}$ is
of type $A_2A_2$ (with $N_{e'}=3$) or $A_3A_1$ (with $N_{e'}$ a power of $2$) so that 
$\cs_1(W_{e'})=\cs(W_{e'})$. Moreover, if $e'\in V^*$ is isolated and $W_{e'}$ is of type $A_2A_2$ then (i) 
holds for this $e'$. We see that $\cs^r_2(W)=\cs_1(W)$ if $r\ne2$ and $\cs^2_2(W)=\cs_2(W)$.  

We prove (f). If $e'\in V^*$ is isolated and $W_{e'}$ is not of type $G_2$ (with $N_{e'}=1$) then from 1.3 we
see that $W_{e'}$ is of type $A_2$ (with $N_{e'}=3$) or $A_1A_1$ (when $N_{e'}=2$) so that 
$\cs_1(W_{e'})=\cs(W_{e'})$. Moreover, if $e'\in V^*$ is isolated and $W_{e'}$ is of type $A_1A_1$, then (i)
holds for this $e'$. We see that $\cs^r_2(W)=\cs_1(W)$ if $r\ne3$ and $\cs^3_2(W)=\cs_2(W)$.  

This proves (a)-(h). From (a)-(h) we deduce:

(j) We have $\cs_2(W)=\cs_2^2(W)\cup\cs_2^3(W)$. If $r\in\cp-\{2,3\}$ then $\cs_2^r(W)=\cs_1(W)$.

The following result can be verified by computation.

(k) If $R$ is of type $E_7$ then $\cs^2_2(W)-\cs_1(W)=\{84_{15}\}$. If $R$ is of type $E_8$ then 
$\cs^2_2(W)-\cs_1(W)=\{1050_{10},840_{14},168_{24},972_{32}\}$ and $\cs^3_2(W)-\cs_1(W)=\{175_{12}\}$. If 
$R$ is of type $F_4$ then $\cs^2_2(W)-\cs_1(W)=\{9_6,4_7,4_8,2_{16}\}$. If $R$ is of type $G_2$ then 
$\cs^3_2(W)-\cs_1(W)=\{1_3\}$.
\nl
(In each case we specify a representation $E$ by a symbol $d_n$ where $d$ is the degree of $E$ and $n=n_E$.
For type $F_4$ and $G_2$ the specified representations are uniquely determined by the additional condition 
that they are not in $\cs_1(W)$.)

(l) $\cs_2^2(W)\cap\cs_2^3(W)=\cs_1(W)$. 
\nl
The inclusion $\cs_1(W)\sub\cs_2^2(W)\cap\cs_2^3(W)$ is obvious. The reverse inclusion for $R$ of type 
$\ne E_8$ follows from the fact that for such $R$ we have either $\cs_2^2(W)=\cs_1(W)$ or 
$\cs_2^3(W)=\cs_1(W)$, see (a)-(h). Thus we can assume that $R$ is of type $E_8$. In this case the result 
follows from (k).

\subhead 1.6\endsubhead
Let $r\in\cp$. Let $V^*_r=\{e'\in V^*;N_{e'}/r\n\ZZ\}$. Let $\tcs^r_2(W)$ be the 
set of all $E_0\in\Irr(W)$ such that for some $e'\in V^*_r$ and some $E\in\cs^r_2(W_{e'})$ we have 
$E_0=j_{W_{e'}}^W(E)$. (Note that any $E\in\cs^r_2(W_{e'})$ is good.) Note that
$\cs^r_2(W)\sub\tcs^r_2(W)$ (take $e'=0$ in the definition of $\tcs^r_2(W)$). We show:

(a) $\cs_2(W)\sub\tcs^r_2(W)$.
\nl
We can assume that $(V,V^*,R,\che R)$ is irreducible. Let $E_0\in\cs_2(W)$. We must show that 
$E_0\in\tcs^r_2(W)$. By 1.1(e) we can find an isolated $e'\n V^*$ and $\tE\in\cs_1(W_{e'})$ such that
$E_0=j_{W_{e'}}^W(\tE)$. If $N_{e'}/r\n\ZZ$ then we have $E_0\in\tcs^r_2(W)$ since
$\cs_1(W_{e'})\sub\cs^r_2(W_{e'})$. If $N_{e'}$ is a power of $r$ then from definitions we have 
$E_0\in\cs^r_2(W)$ hence $E_0\in\tcs^r_2(W)$. Thus we may assume that $N_{e'}$ is not a power of $r$ and is
$N_{e'}/r\in\ZZ$. This forces $R$ to be of type $E_8$ and $W_{e'}$ to be of type $A_5A_2A_1$ (see 1.3); we
then have $N_{e'}=6$ and $r\in\{2,3\}$. In particular we must have $\tE\in\cs(W_{e'})$. If $\tE$ is not the 
sign representation of $W_{e'}$ then we have $\tE=j_{W_{e'_1}}^{W_{e'}}(\text{sign})$ for some $e'_1\in V^*$
such that $W_{e'_1}$ is a proper parabolic subgroup of $W_{e'}$. Replacing $W_{e'_1}$ by a $W$-conjugate we 
can assume that $W_{e'_1}$ is a proper parabolic subgroup of $W$ so that 
$j_{W_{e'}}^W(\text{sign})\in\cs(W)$ and in particular, $E_0\in\tcs^r_2(W)$.
Thus we can assume that $\tE$ is the sign representation of $W_{e'}$. We have $W_{e'}\sub W_{e'_2}$ where 
$W_{e'_2}$ is of type $E_7A_1$ and by the definition of $\cs_1(W_{e'_2})$ we have 
$$\tE_2:=j_{W_{e'}}^{W_{e'_2}}(\text{sign})\in\cs_1(W_{e'_2}).$$
If $r=3$ we have $e'_2\in V^*_r$ hence $E_0=j_{W_{e'_2}}^W(\tE_2)\in\tcs^r_2(W)$. We have 
$W_{e'}\sub W_{e'_3}$ where $W_{e'_3}$ is of type $E_6A_2$ and by the definition of $\cs_1(W_{e'_3})$ we 
have $\tE_3:=j_{W_{e'}}^{W_{e'_3}}(\text{sign})\in\cs_1(W_{e'_3})$. If $r=2$ we have $e'_3\in V^*_r$ hence 
$E_0=j_{W_{e'_3}}^W(\tE_3)\in\tcs^r_2(W)$. This completes the proof of (a).

We show:

(b) $\tcs^r_2(W)\sub\cs_2(W)$.
\nl
We can assume that $(V,V^*,R,\che R)$ is irreducible. Let $E_0\in\tcs^r_2(W)$. We must show that 
$E_0\in\cs_2(W)$. Assume first that $r\n\{2,3\}$. Then by results in 1.5 we have $E\in\cs_1(W_{e'})$ hence 
by 1.1(c) we have $E_0\in\cs_2(W)$. Next we assume that $r=3$. If $W_{e'}\ne W$ then by results in 1.5 we 
have $E\in\cs_1(W_{e'})$ hence by 1.1(c) we have $E_0\in\cs_2(W)$. Thus we can assume that $W_{e'}=W$ so 
that $E_0=E\in\cs^r_2(W)$. Since $\cs^r_2(W)\sub\cs_2(W)$ we see that $E_0\in\cs_2(W)$.

We now assume that $r=2$. We can find $e'\in V^*_r$ and $E\in\cs^r_2(W_{e'})$ such that 
$E_0=j_{W_{e'}}^W(E)$. 
We can find an isolated $e'_1\in V^*$ such that $N_{e'_1}$ is odd, $R_{e'}\sub R_{e'_1}$ and $R_{e'}$ is 
rationally closed in $R_{e'_1}$. Let $E'=j_{W_{e'}}^{W_{e'_1}}(E)$. Since $E\in\cs_2(W_{e'})$ we have
$E'\in\cs_2(W_{e'_1})$, see 1.1(g) and $E_0=j_{W_{e'_1}}^W(E')$. It is then enough to prove the following 
statement:

(c) If $e'\in V^*_r$ is isolated $(r=2)$ and $E\in\cs_2(W_{e'})$ then $E_0=j_{W_{e'}}^W(E)\in\cs_2(W)$. 
\nl
If $W_{e'}=W$ then $E_0=E\in\cs_2(W)$, as required. If $R$ is of type $A_n,B_n,C_n,D_n$ then in (c) we have
automatically $W_{e'}=W$ hence (c) holds in these cases. Thus we can assume in (c) that $R$ is of exceptional
type and $W_{e'}\ne W$. Then $W_{e'}$ is of the following type: 
$A_2A_2A_2$ (if $R$ is of type $E_6$); $A_5A_2$ (if $R$ is of type $E_7$); $A_4A_4$ or $A_8$ or $E_6A_2$ (if
$R$ is of type $E_8$); $A_2A_2$, as in 1.3(h) (if $R$ is of type $F_4$); $A_2$, as in 1.3(i) (if $R$ is of
type $G_2$). In each case we have $\cs_2(W_{e'})=\cs_1(W_{e'})$, see 1.5. Thus $E\in\cs_1(W_{e'})$. Using
1.1(e) we see that $E_0\in\cs_2(W)$. This proves (c) hence (b).

\mpb

Combining (a),(b) we obtain

(d) $\tcs^r_2(W)=\cs_2(W)$.
\nl
In the case where $r=0$, we set $V^*_0=V^*$, $\cs^0_2(W)=\cs_1(W)$, $\tcs^0_2(W)=\cs_2(W)$.

\head 2. The strata of $G$\endhead
\subhead 2.1\endsubhead 
We return to the setup of 0.1. Thus $G$ is a connected reductive algebraic group over $\kk$. Let $\ct$ be 
``the'' maximal torus of $G$; let $X=\Hom(\ct,\kk^*)$, $Y=\Hom(\kk^*,\ct)$, $V=\QQ\ot X$, $V^*=\QQ\ot Y$. 
We have an obvious perfect bilinear pairing $\la,\ra:V\T V^*@>>>\QQ$.
Let $R\sub V$ be the set of roots and let $\che R\sub V^*$ be the set of corrots. Then $(V,V^*,R,\che R)$
is as in 1.1. The associated Weyl group $W$ (as in 1.1) that is, the Weyl group of $G$, can be viewed as an
indexing set for the orbits of $G$ acting diagonally on $\cb\T\cb$; we denote by 
$\co_w$ the orbit corresponding to $w\in W$. Note that $W$ is naturally a Coxeter group.

Let $g\in G$. Let $W_g$ be the Weyl group of the connected reductive group $H:=Z_G(g_s)^0$. We can view 
$W_g$ as a subgroup of $W$ as follows. Let $\b$ be a Borel subgroup of $H$ and let $T$ be a maximal torus of
$\b$. We define an isomorphism $b_{T,\b}:N_HT/T@>\si>>W_g$ by $n'T\m H\text{-orbit of }(\b,n'\b n'{}\i)$. 
Similarly for any $B\in\cb$ such that $T\sub B$ we define an isomorphism $a_{T,B}:N_GT/T@>\si>>W$ by 
$n'T\m G\text{-orbit of }(B,n'Bn'{}\i)$. Now assume that
$B\in\cb$ is such that $B\cap H=\b$. We define an imbedding $c_{T,\b,B}:W_g@>>>W$ as the 
composition $W_g@>b\i_{T,\b}>>N_HT/T@>>>N_GT/T@>a_{T,B}>>W$ where the middle map is the obvious imbedding. 
If $B'\in\cb$ also satisfies $B'\cap H=\b$ then we have $B'=nBn\i$ for some $n\in N_GT$ and from the 
definitions we have $c_{T,\b,B'}(w)=a_{T,B}(nT)c_{T,\b,B}(w)a_{B,T}(nT)\i$ for any $w\in W_g$. Thus 
$c_{T,\b,B}$ depends (up to composition with an inner automorphism of $W$) only on $T,\b$ and we can denote 
it by $c_{T,\b}$. Since the set of pairs $T,\b$ as above form a homogeneous space for the connected group 
$H$ we see that $c_{T,\b}$ is independent of $T,\b$ (up to composition with an inner automorphism of $W$) 
hence it does not depend on any choice. We see that there is a well defined collection $\cc$ of imbeddings 
$W_g@>>>W$ so that any two of them differ only by composition by an inner automorphism of $W$. 

Define $\r\in\Irr(W_g)$ by the condition that under the Springer correspondence for $H$, $\r$ corresponds to
the $H$-conjugacy class of $g_u$ and the trivial local system on it. We choose $f\in\cc$; then we can view 
$\r$ as an irreducible representation of $f(W_g)$, a subgroup of $W$ such that $f(W_g)=W_{e'}$ for some 
$e'\in V^*_p$, see 1.6. By \cite{\UNI, 1.4} we have $\r\in\cs^p_2(f(W_g))$, see 1.5, 1.6. Hence 
$\ti\r:=j_{f(W_g)}^W(\r)\in\tcs^p_2(W)$ is well defined. Since $\tcs^p_2(W)=\cs_2(W)$, see 1.6, we have 
$\ti\r\in\cs_2(W)$. This is independent of the choice of $f$ since $f$ is well defined up to composition by 
an inner automorphism of $W$.

\subhead 2.2\endsubhead
Let $g\in G$. Let $d=d_g=\dim\cb_g$. The imbedding $h_g:\cb_g@>>>\cb$ induces a linear map
$h_{g*}:H_{2d}(\cb_g)@>>>H_{2d}(\cb)$. Now $H^{2d}(\cb_g),H^{2d}(\cb)$ carry natural $W$-actions, see 
\cite{\ICC}, and this induces natural $W$-actions on $H_{2d}(\cb_g),H_{2d}(\cb)$ which are compatible with 
$h_{g*}$. Hence $W$ acts naturally on the subspace $h_{g*}(H_{2d}(\cb_g))$ of $H_{2d}(\cb)$.

The following result gives an alternative description of the map $g\m\ti\r$ (in 2.1) from $G$ to $\Irr W$. 

(a) {\it The $W$-submodule $h_{g*}(H_{2d}(\cb_g))$ of $H_{2d}(\cb)$ is isomorphic to the $W$-module 
$\QQ_l\ot\ti\r$ where $\r,\ti\r$ are associated to $g$ as in 2.1.}
\nl
First we note that $h_{g*}(H_{2d}(\cb_g))\ne0$; indeed it is clear that for any irreducible component $D$ of
$\cb_g$ (necesarily of dimension $d$) the image of the fundamental class of $D$ under $h_{g*}$ is nonzero 
(we ignore Tate twists). Let $\cb'$ be the variety of Borel subgroups of $Z_G(g_s)^0$. Let 
$\cb'_{g_u}=\{\b\in\cb';g_u\in\b\}$. Then $\dim\cb'=d$ and $W_g$ (see 2.1) acts naturally on 
$H_{2d}(\cb'_{g_u})$; from the definitions, the $W$-module $H_{2d}(\cb_g)$ is isomorphic to 
$\Ind_{W_g}^WH_{2d}(\cb'_{g_u})$. From the definitions we have $n_\r=d$ and the $W_g$-module 
$H_{2d}(\cb'_{g_u})$ is of the form $\op_{i\in[1,s]}(\bbq\ot E_i)^{\op c_i}$ where $E_i\in\Irr(W_g)$, 
$c_i\in\NN$ satisfy $E_1=\r$, $c_1=1$ and $n_{E_i}>d$ for $i>1$. It follows that the $W$-module 
$H_{2d}(\cb_g)$ is of the form $\op_{i\in[1,s]}(\Ind_{W_g}^W(\bbq\ot E_i))^{\op c_i}$. Now 
$\Ind_{W_g}^W(\bbq\ot E_1)$ contains $\bbq\ot\ti\r$ with multiplicity $1$ and all its other irreducible 
constituents are of the form $\bbq\ot E$ with $n_E>d$; moreover, for $i>1$, any irreducible constituent $E$ 
of $\Ind_{W_g}^W(\bbq\ot E_i)$ satisfies $n_E>d$. Thus the $W$-module $H_{2d}(\cb_g)$ contains 
$\bbq\ot\ti\r$ with multiplicity $1$ and all its other irreducible constituents are of the form $\bbq\ot E$ 
with $n_E>d$; these other irreducible constituents are necessarily mapped to $0$ by $h_{g*}$ and the 
irreducible constituent isomorphic to $\bbq\ot\ti\r$ is mapped injectively by $h_{g*}$ since $h_{g*}\ne0$. 
It follows that the image of $h_{g*}$ is isomorphic to $\bbq\ot\ti\r$ as a $W$-module. This proves (a).

\subhead 2.3\endsubhead
By 2.1, 2.2 we have a well defined map $\ph:G@>>>\cs_2(W)$, $g\m\ti\r$ where 
$\QQ_l\ot\ti\r=h_{g*}((H_{2d_g}(\cb_g))$ (notation of
2.1, 2.2). The fibres $G_E=\ph\i(E)$ of $\ph$ ($E\in\cs_2(W)$) are called the {\it strata} of $G$. They are 
clearly unions of conjugacy 
classes of $G$. Note the strata of $G$ are indexed by the finite set $\cs_2(W)$ which depends only on the 
Weyl group $W$ and not on the underlying root system (see 1.1(b)) or on the characteristic of $\kk$. 

One can show that any stratum of $G$ is a union of pieces in the partition of $G$ defined in 
\cite{\ICC, 3.1}; in particular, it is a constructible subset of $G$. 

\subhead 2.4\endsubhead 
We have the following result.

(a) {\it Any stratum $G_E$ ($E\in\cs_2(W)$) of $G$ is a (non-empty) union of $G$-conjugacy classes of fixed 
dimension, namely $2\dim\cb-2n$ where $n=n_E$, see 0.2. At most one $G$-conjugacy class in $G_E$ is 
unipotent.}
\nl
Since $\cs_2(W)=\tcs^p_2(W)$, see 1.6, we have
$E\in\tcs^p_2(W)$. Hence there exists $e'\in V^*_p$ and $\r\in\cs^p_2(W_{e'})$ such 
that $E=j_{W_{e'}}^W(\r)$. We can find a semisimple element of finite order $s\in G$ such that $W_s$ (viewed
as a subgroup of $W$ as in 2.1) is equal to $W_{e'}$. By \cite{\UNI, 1.4} we can find a unipotent element 
$u$ in $Z_G(s)^0$ such that $\r$ is the Springer representation of $W_s$ defined by $u$ and the trivial
local system on its $Z_G(s)^0$-conjugacy class. Then $E=\ph(su)$ so that $G_E\ne\emp$.
Let $\g$ be a $G$-conjugacy class in $G_E$. Let $g\in\g$. Let $\r$ (resp. $\ti\r$) be the irreducible
representation of $W_g$ (resp. $W$) defined by $g_u$ as in 2.1. Let $n_\r,n_{\ti\r}$ be as in 0.2. By the 
definition of $\ti\r$ we have $n_\r=n_{\ti\r}$. By assumption we have $\ti\r=E$ hence $n_{\ti\r}=n$ and 
$n_\r=n$. By a known property of Springer's representations, $n_\r$ is equal to the dimension of the variety
of Borel subgroups of $Z_G(g_s)^0$ that contain $g_u$; hence by a result of Steinberg (for $p=0$) and 
Spaltenstein \cite{\SPA, 10.15} (for any $p$), $n_\r$ is equal to 
$$(\dim(Z_{Z_G(g_s)^0}(g_u)^0-\rk(Z_G(g_s)^0))/2=(\dim(Z_G(g)^0)-\rk(G))/2.$$
It follows that $(\dim(Z_G(g)^0)-\rk(G))/2=n$ and the desired formula for $\dim\g$ follows. Now assume that 
$\g,\g'$ are two unipotent $G$-conjugacy classes contained in $G_E$. Then the Springer representation of $W$
associated to $\g$ is the same as that associated to $\g'$, namely $E$. By properties of Springer 
representations, it follows that $\g=\g'$. This proves (a).

\subhead 2.5\endsubhead 
In this and the next subsection we assume that $W$ is irreducibble.
Let $r\in\cp\cup\{0\}$. Let $G^r$ be a connected reductive group of the same type as $G$
over an algebraically closed field of characteristic $r$, whose Weyl group is identified with $W$. 
Let $\cu^r$ be the set of unipotent classes of $G^r$. By \cite{\UNI, 1.4} we have a canonical bijection
$$\ps^r:\cu^r@>\si>>\cs^r_2(W)$$
which to a unipotent class $\g$ associates the Springer representation of $W$ corresponding to $\g$ and the
constant local system on $\g$. We define an imbedding $h^r:\cu^0@>>>\cu^r$ as the composition 
$$\cu^0@>\ps^0>>\cs^0_2(W)=\cs_1(W)@>>>\cs^r_2(W)@>(\ps^r)\i>>\cu^r$$
where the unnamed map is the inclusion.

Consider the relation $\cong$ on $\sqc_{r\in\cp}\cu^r$ for which $x\in\cu^r$, $y\in\cu^{r'}$ (where 
$r,r'\in\cp$) satisfy $x\cong y$ if either $r=r'$ and $x=y$ or $r\ne r'$ and $x=h^r(z)$, $y=h^{r'}(z)$ 
for some $z\in\cu^0$. We show that $\cong$ is an equivalence relation.
It is enough to show that if $x\in\cu^r$, $y\in\cu^{r'}$, $u\in\cu^{r''}$ are such that
$r\ne r'$, $r'\ne r''$  and $x=h^r(z)$, $y=h^{r'}(z)$,  $y=h^{r'}(\tz)$, $u=h^{r''}(\tz)$ for some
$z\in\cu^0,\tz\in\cu^0$ then $x\cong u$. From $h^{r'}(z)=h^{r'}(\tz)$ and the injectivity of $h^{r'}$ we
have $z=\tz$. Thus, if $r\ne r''$ we have $x\cong u$ while if $r=r''$ we have $x=u$. Thus, $\cong$ is indeed
an equivalence relation.

Let $\cu^*$ be $\sqc_{r\in\cp}\cu^r$ modulo the equivalence relation $\cong$.
Let $\sqc_{r\in\cp}\cu^r@>>>\cs_2(W)$ be the map whose restriction to $\cu^r$ is $\ps^r$ followed by 
the inclusion $\cs^r_2(W)\sub\cs_2(W)$ (for any $r$). We show:

(a) {\it This map induces a bijection $\ps^*:\cu^*@>\si>>\cs_2(W)$.}
\nl
To show that $\ps^*$ is a well defined map it is enough to verify that if $z\in\cu^0$ then for any
$r,r'\in\cp$ we have 
$\ps^rh^r(z)=\ps^{r'}h^{r'}(z)$ in $\cs_2(W)$; but both sides of the equality to be verified are equal to
$\ps^0(z)$. Let $E\in\cs_2(W)$. By 1.5(j) there
exists $r\in\cp$ such that $E\in\cs^r_2(W)$ hence $E=\ps^r(x)$ for some $x\in\cu^r$. It follows that 
$\ps^*$ is surjective. We show that $\ps^*$ is injective. It is enough to show that 

(b) if $x\in\cu^r$, $y\in\cu^{r'}$ ($r,r'\in\cp$ distinct) satisfy $\ps^r(x)=\ps^{r'}(y)$ then there exists 
$z\in\cu^0$ such that $x=h^r(z)$, $y=h^{r'}(z)$.
\nl
If $r\ne\{2,3\}$ then $\cs^r_2(W)=\cs_1(W)$ hence $\ps^r(x)=\ps^0(z)$ for some $z\in\cu^0$. We then have
$\ps^{r'}(y)=\ps^0(z)$. It follows that $h^r(z)=x$, $h^{r'}(z)=y$, as required. Similarly, if $r'\ne\{2,3\}$
then the conclusion of (b) holds. Thus we can assume that $r\in\{2,3\}$, $r'\in\{2,3\}$. Since $r\ne r'$ we
have $\{r,r'\}=\{2,3\}$. Hence $\ps^r(x)=\ps^{r'}(y)\in\cs^2_2(W)\cap\cs^3_2(W)=\cs_1(W)$; the last equality
follows from 1.5(l). Thus we have 
$\ps^r(x)=\ps^{r'}(y)=\ps^0(z)$ for some $z\in\cu^0$. It follows that $h^r(z)=x$, $h^{r'}(z)=y$, as required.

\mpb

From (a) we deduce:

(c) {\it The strata of $G$ are naturally indexed by the set $\cu^*$.}
\nl
The proof of (a) shows also that 
$\cu^*$ is equal to $\cu^2\sqc\cu^3$ with the identification of $h^2(z),h^3(z)$ for any $z\in\cu^0$.

We can now state the following result.

(d) {\it Let $E\in\cs_2(W)$. Then for some $r\in\cp$, the stratum $G^r_E$ contains a unipotent class. In 
fact, $r$ can be assumed to be $2$ or $3$.}
\nl
Under (a), $E$ corresponds to an element of $\cu^*$ which is the equivalence class of some element 
$\g\in\cu^r$ with $r\in\{2,3\}$. Let $g\in G^r$ be an element in the unipotent conjugacy class $\g$. From
the definitions we see that $g\in G^r_E$. This proves (d).

\subhead 2.6\endsubhead 
We show that the set $\cu^*$ has a natural partial order.
If $\cs^r_2(W)=\cs_1(W)$ (type $A$ and $E_6$) we have $\cu^*=\cu^0$
which has a natural partial order defined by the closure relation of unipotent classes in $G^0$. If 
$\cs_1(W)\ne\cs^r_2(W)$ for a unique $r\in\cp$ (type $\ne A,E_6,E_8$), we have $\cu^*=\cu^r$ which has a 
natural partial order defined by the closure relation of unipotent classes in $G^r$.
Assume now that $G$ is of type $E_8$. Then we can identify $\cu^2,\cu^3$ with subsets of $\cu^*$ whose union
is $\cu^*$ and whose intersection is $\cu^0$. Both subsets $\cu^2,\cu^3$ have natural partial orders defined
by the closure relation of unipotent classes in $G^2$ and $G^3$. If $\g,\g'\in\cu^*$ we say that $\g\le\g'$ 
if there exists a sequence $\g=\g_0,\g_1,\do,\g_s=\g'$ in $\cu^*$ such that for any $i\in[1,s]$ there exists
$r\in\{2,3\}$ such that 

(a) $\g_{i-1}\in\cu^r,\g_i\in\cu^r$, $\g_{i-1}\le\g_i$ in the partial order of unipotent classes in $G^r$;
\nl
note that if for some $i$, (a) holds for both $r=2$ and $r=3$ then we have 
$\g_{i-1}\in\cu^0,\g_i\in\cu^0$, $\g_{i-1}\le\g_i$ in the partial order of unipotent classes in $G^0$. One 
can show that this partial order on $\cu^*$ induces the usual partial orders on the subsets $\cu^2,\cu^3$,
$\cu^0$.

\subhead 2.7\endsubhead 
Let $W_a$ be the semidirect product of $W$ with the subgroup of $V$ generated by $R$ (an affine Weyl group);
let ${}'W_a$ be the semidirect product of $W$ with the subgroup of $V^*$ generated by $\che R$ (another affine
Weyl group). We consider four triples: 

(a) $(\cs(W), X_0, Z_0)$

(b) $(\cs_1(W), X_1, Z_1)$

(c) $({}'\cs_1(W), {}'X_1, {}'Z_1)$

(d) $(\cs_2(W), X_2, Z_2)$
\nl
where $X_0,X_1,{}'X_1$ is the set of two-sided cells in $W,W_a,{}'W_a$ respectively, $Z_0$ is the set of 
special unipotent classes in $G$ with $p=0$, $Z_1$ is the set of unipotent classes in $G$ with $p=0$,
${}'Z_1$ is the set of unipotent classes in the Langlands dual $G^*$ of $G$ with $p=0$, $Z_2$ is the set of
strata of $G$ with $p=0$ and $X_2$ remains to be defined.
The three sets in each of these four triples are in canonical bijection with each other (assuming that $X_2$
has been defined). Moreover, each set in (a) is naturally contained in the corresponding set in (b) and 
(replacing $G$ by $G^*$) in the corresponding set in (c); each set in (b) is contained in the corresponding 
set in (d) and (replacing $G$ by $G^*$) each set in (c) is contained in the corresponding set in (d). It 
remains to define $X_2$. It seems plausible that the (trigonometric) 
double affine Hecke algebra $\HH$ associated by Cherednik to $W$ has a natural filtration by two-sided 
ideals whose successive subquotients can be called two-sided cells and form the desired set $X_2$.
The inclusion of the Hecke algebra of $W_a$ and that of ${}'W_a$ into $\HH$ should induce the imbeddings 
$X_1\sub X_2$, ${}'X_1\sub X_2$ and $X_2$ should be in natural bijection with
$\cs_2(W)$ and with the set of strata of $G$.

\head 3. Examples\endhead
\subhead 3.1\endsubhead
We write the adjoint group of $G$ as a product $\prod_iG_i$ where each $G_i$ is simple with Weyl group $W_i$
so that $W=\prod_iW_i$. Let $E\in\cs_2(W)$. We have $E=\bxt_iE_i$ where $E_i\in\cs_2(W_i)$. Now $G_E$ is the
the inverse image of $\prod_i(G_i)_{E_i}$ under the obvious map $G@>>>\prod_iG_i$.

When $E$ is the sign representation of $W$ then $G_E$ is the centre of $G$; when $E$ is the unit
representation of $W$, $G_E$ is the set of elements of $G$ which are regular in the sense of Steinberg
\cite{\ST}.

By 2.5(a) and 2.6 applied to $G_i$, the set $\cs_2(W_i)$ has a natural partial order. Since
$\cs_2(W)$ can be identified as above with $\prod_i\cs_2(W_i)$, $\cs_2(W)$ is naturally a partially
ordered set (a product of partially ordered sets). Hence by 2.3 the set of strata of $G$ is 
naturally a partially ordered set.

\subhead 3.2\endsubhead 
Assume that $G=GL(V)$ where $V$ is a $\kk$-vector space of dimension $n\ge1$. Let $g\in G$. For any 
$x\in\kk^*$ let $V_x$ be the generalized $x$-eigenspace of $g:V@>>>V$ and let 
$\l^x_1\ge\l^x_2\ge\l^x_3\ge\do$ be the sequence in $\NN$ whose nonzero terms are the sizes of the Jordan 
blocks of $x\i g:V_x@>>>V_x$. Let ${}^g\l$ be the sequence ${}^g\l_1\ge{}^g\l_2\ge{}^g\l_3\ge\do$ given by 
${}^g\l_j=\sum_{x\in\kk^*}\l^x_j$. Now $g\m{}^g\l$ defines a map from $G$ onto the set of partitions of $n$.
From the definitions we see that the fibres of this map are exactly the strata of $G$. If $g\in G$ and 
${}^g\l=(\l_1,\l_2,\l_3,\do)$ then 
$$\dim(\cb_g)=\sum_{k\ge1}(n-(\l_1+\l_2+\do+\l_k)).$$

\subhead 3.3\endsubhead 
Repeating the definition of sheets in a semisimple Lie algebra over $\CC$ (see \cite{\BO}) one can define
the sheets of $G$ as the maximal irreducible subsets of $G$ which are unions of conjugacy classes of fixed
dimension. One can show that if $G$ is as in 4.2, the sheets of $G$ are the same as the strata of $G$, as 
described in 4.2. (In this case, the sheets of $G$, or rather their Lie algebra analogue, are described in 
\cite{\PE}. They are smooth varieties.) This is not true for a general $G$ (the sheets of $G$ do not usually
form a partition of $G$; the strata of $G$ are not always irreducible). In \cite{\CAR} it is shown that if 
$p$ is $0$ or a good prime for $G$ then any stratum is a union of sheets and that the closure of a stratum is
not necessarily a union of strata, even if $G$ is of type $A$. 

\subhead 3.4\endsubhead 
In the next few subsections we will describe explicitly the strata of $G$ when $G$ is a symplectic or
special ortogonal group.

Given a partition $\nu=(\nu_1\ge\nu_2\ge\do)$, a {\it string} of $\nu$ is a maximal subsequence
$\nu_i,\nu_{i+1},\do,\nu_j$  of $\nu$ consisting of equal $>0$ numbers; the string is said to have an odd
origin if $i$ is odd and an even origin if $i$ is even.

For an even $N\in\NN$ let $Z^1_N$ be the set of partitions $\nu=(\nu_1\ge\nu_2\ge\do)$ of $N$ such that any 
odd number appears an even number of times in $\nu$. We show:

(a) {\it There is a canonical bijection $Z^1_N\lra BP^{N/2}_{1,1}$ (notation of 0.2).}
\nl
To $\nu\in Z^1_N$ we associate $\l=(\l_1,\l_2,\l_3,\do)$ as follows: each string $2a,2a,\do,2a$ in $\nu$ is 
replaced by $a,a,\do,a$ of the same length; each string $2a+1,2a+1,\do,2a+1$ (necessarily of even length) in
$\nu$ is replaced by $a,a+1,a,a+1,\do,a,a+1$ of the same length. The resulting entries form a bipartition 
$\l\in BP^{N/2}_{1,1}$. Now $\nu\m\l$ establishes the bijection (a).

\mpb

For an even $N\in\NN$ let $Z^2_N$ be the set of partitions $\nu=(\nu_1\ge\nu_2\ge\do)$ of $N$ such that any
odd number appears an even number of times in $\nu$ and any even $>0$ number which appears an even $>0$ 
number of times in $\nu$ has an associated label $0$ or $1$. We show:

(b) {\it There is a canonical bijection $Z^2_N\lra BP^{N/2}_{2,2}$ (notation of 0.2).}
\nl
To $\nu\in Z^2_N$ we associate $\l=(\l_1,\l_2,\l_3,\do)$ as follows: each string $2a,2a,\do,2a$ of odd 
length or of even length and label $1$ in $\nu$ is replaced by $a,a,\do,a$ of the same length; each string 
$2a,2a,\do,2a$ of even length and label $0$ in $\nu$ is replaced by $a-1,a+1,a-1,a+1,\do,a-1,a+1$ of the 
same length; each string $2a+1,2a+1,\do,2a+1$ (necessarily of even length) in $\nu$ is replaced by 
$a,a+1,a,a+1,\do,a,a+1$ of the same length. The resulting entries form a bipartition $\l\in BP^{N/2}_{2,2}$.
Now $\nu\m\l$ establishes the bijection (b).

Assume for example that $N=6$. The bijection (b) is:

$(6\do)\lra(3\do)$

$(42\do)\lra(21\do)$

$(411\do)\lra(201\do)$

$(33\do)\lra(12\do)$

$(222\do)\lra(111\do)$

$((22)_111\do)\lra(1101\do)$

$((22)_0110\do)\lra(0201\do)$

$(21111\do)\lra(10101\do)$

$(111111\do)\lra(010101\do)$.
\nl
Here we write $\do$ instead of $000\do$. (Compare \cite{\LSS, 6.1}.)

\subhead 3.5\endsubhead 
Assume that $G=Sp(V)$ where $V$ is a $\kk$-vector space of dimension $N$ with a fixed nondegenerate
symplectic form. 

Let $g\in G$. For any $x\in\kk^*$ let $V_x$ be the generalized $x$-eigenspace of $g:V@>>>V$. Let
$d_x=\dim V_x$. For any $x\in\kk^*$ such that $x^2\ne1$ let $\l^x_1\ge\l^x_2\ge\l_3^x\ge\do$ be the 
partition of $d_x$ whose nonzero terms are the sizes of the Jordan blocks of $x\i g:V_x@>>>V_x$. 

For $x\in\kk^*$ such that $x^2=1$ let $\nu^x\in Z^1_{d_x}$ (if $p\ne2$) and
$\nu^x\in Z^2_{d_x}$ (if $p=2$) be again the partition of $d_x$ 
whose nonzero terms are the sizes of the Jordan blocks of the unipotent element $x\i g\in Sp(V_x)$. (When 
$p=2$, $\nu^x$ should also include a labelling with $0$ and $1$ associated to $x\i g\in Sp(V_x)$ as in 
\cite{\WEUII, 1.4}.) Let $\l^x=(\l^x_1,\l^x_2,\l_3^x,\do)$
 be the bipartition of $d_x/2$ associated to $\nu^x$ by 3.4(a),(b). Thus
$\l^x\in BP^{d_x/2}_{1,1}$ (if $p\ne2$), $\l^x\in BP^{d_x/2}_{2,2}$ (if $p=2$).
Note that $\l^x$ is the bipartition such that the Springer representation attached to the unipotent element 
$x\i g\in Sp(V_x)$ (an irreducible representation of the Weyl group of type $B_{d_x/2}$) is indexed in the 
standard way by $\l^x$.
Define ${}^g\l=({}^g\l_1,{}^g\l_2,{}^g\l_3,\do)$ by ${}^g\l_j=\sum_x\l^x_j$ where $x$ runs over
a set of representatives for the orbits of the involution $a\m a\i$ of $\kk^*$. Note that 
${}^g\l\in BP^{N/2}_{2,2}$. Thus we have defined a (surjective) map $g\m{}^g\l$, $G@>>>BP^{N/2}_{2,2}$.
From the definitions we see that the fibres of this map are exactly the strata of $G$. 

If $g\in G$ and ${}^g\l=(\l_1,\l_2,\l_3,\do)$ then 
$$\dim(\cb_g)=\sum_{k\ge1}((N/2)-(\l_1+\l_2+\do+\l_k)).\tag a$$ 

We now consider the case where $N=4$. In this case we have $\cs_2(W)=\Irr(W)$; hence there are five strata. 
One stratum is the union of all conjugacy classes of dimension $8$ (it corresponds to the unit 
representation); one stratum is the union of all conjugacy classes of dimension $6$ (it corresponds to the 
reflection representation of $W$). There are two strata which are unions of conjugacy classes of dimension 
$4$ (they correspond to the two one dimensional representations of $W$ other than unit and sign); if $p=2$
both these strata are single unipotent classes; if $p\ne2$ one of these strata is a semisimple class and the
other is a unipotent class times the centre of $G$. The centre of $G$ is a stratum (it corresponds to the 
sign representation of $W$).

The results in this subsection show that under the standard identification $\Irr(W)=BP^{N/2}$, we have
$$\cs_2(W)=BP^{N/2}_{2,2}.\tag b$$
Under this identification the map $g\m{}^g\l$, $G@>>>BP^{N/2}_{2,2}$ becomes the map $g\m E$ where
$g\in G_E$.

\subhead 3.6\endsubhead 
For $N\in\NN$ let ${}'Z^1_N$ be the set of partitions $\nu=(\nu_1\ge\nu_2\ge...)$ such that any even $>0$
number appears an even number of times in $\nu$ and $\nu_1+\nu_2+\do=N$.

(a) {\it If $N$ is odd, then there is a canonical bijection ${}'Z^1_N\lra BP^{(N-1)/2}_{2,0}$.}
\nl
To $\nu\in{}'Z^1_N$ we associate $\l=(\l_1,\l_2,\l_3,\do)$ as follows: each string $2a,2a,\do,2a$ of $\nu$ 
(necessarily of even length) is replaced by $a-1,a+1,a-1,a+1,\do,a-1,a+1$ of the same length (if the string 
has odd origin) or by $a,a,\do,a$ of the same length (if the string has even origin); each string 
$2a+1,2a+1,\do,2a+1$ of $\nu$ is replaced by $a,a+1,a,a+1,\do$ of the same length (if the string has odd 
origin) or by $a+1,a,a+1,a,\do$ of the same length (if the string has even origin). The resulting entries 
form a bipartition $\l\in BP^{(N-1)/2}_{2,0}$. Now $\nu\m\l$ establishes the bijection (a).

\mpb

(b) {\it If $N$ is even, then there is a canonical bijection ${}'Z^1_N\lra BP^{N/2}_{0,2}$.}
\nl
To $\nu\in{}'Z^1_N$ we associate $\l=(\l_1,\l_2,\l_3,\do)$ as follows: each string $2a,2a,\do,2a$ of $\nu$ 
(necessarily of even length) is replaced by $a-1,a+1,a-1,a+1,\do,a-1,a+1$ of the same length (if the string 
has even origin) or by $a,a,\do,a$ of the same length (if the string has odd origin); each string 
$2a+1,2a+1,\do,2a+1$ of $\nu$ is replaced by $a,a+1,a,a+1,\do$ of the same length (if the string has even 
origin) or by $a+1,a,a+1,a,\do$ of the same length (if the string has odd origin). The resulting entries 
form a bipartition $\l\in BP^{N/2}_{0,2}$. Now $\nu\m\l$ establishes the bijection (b).

\subhead 3.7\endsubhead 
Assume that $p\ne2$ and that $G=SO(V)$ where $V$ is a $\kk$-vector space of odd dimension $N\ge1$ with a 
fixed nondegenerate quadratic form. 

Let $g\in G$. For any $x\in\kk^*$ let $V_x$ be the generalized $x$-eigenspace of $g:V@>>>V$. Let 
$d_x=\dim V_x$. For any $x\in\kk^*$ such that $x^2\ne1$ let $\l^x_1\ge\l^x_2\ge\l_3^x\ge\do$ be the 
partition of $d_x$ whose nonzero terms are the sizes of the Jordan blocks of $x\i g:V_x@>>>V_x$. 

For $x\in\kk^*$ such that $x^2=1$ let $\nu^x\in{}'Z^1_{d_x}$ be again the partition of $d_x$ 
whose nonzero terms are the sizes of the Jordan blocks of the unipotent element $x\i g\in SO(V_x)$. 
Let $\l^x=(\l^x_1,\l^x_2,\l_3^x,\do)$ be the bipartition of $d_x/2$ associated to $\nu^x$ by 3.6(a) if
$x=1$ and by 3.6(b) if $x=-1$. Thus $\l^x\in BP^{(d_x-1)/2}_{2,0}$ if $x=1$,
$\l^x\in BP^{d_x/2}_{0,2}$ if $x=-1$.
Note that $\l^x$ is the bipartition such that the Springer representation attached to the unipotent element 
$x\i g\in SO(V_x)$ (an irreducible representation of the Weyl group of type $B_{(d_x-1)/2}$, if $x=1$, or of
type $D_{d_x/2}$, if $x=-1$) is indexed by $\l^x$. 
Define ${}^g\l=({}^g\l_1,{}^g\l_2,{}^g\l_3,\do)$ by ${}^g\l_j=\sum_x\l^x_j$ where $x$ runs over
a set of representatives for the orbits of the involution $a\m a\i$ of $\kk^*$. Note 
that ${}^g\l\in BP^{(N-1)/2}_{2,2}$.
Thus we have defined a (surjective) map $g\m{}^g\l$, $G@>>>BP^{(N-1)/2}_{2,2}$.
From the definitions we see that the fibres of this map are exactly the strata of $G$. 
Under the identification $\cs_2(W)=BP^{(N-1)/2}_{2,2}$, see 3.5(b), 
the map $g\m{}^g\l$, $G@>>>BP^{(N-1)/2}_{2,2}$ becomes the map $g\m E$ where $g\in G_E$.

If $g\in G$ and ${}^g\l=(\l_1,\l_2,\l_3,\do)$ then 
$$\dim(\cb_g)=\sum_{k\ge1}((N-1)/2-(\l_1+\l_2+\do+\l_k)).$$ 

\subhead 3.8\endsubhead
Assume that $p=2$ and that $G=SO(V)$ where $V$ is a $\kk$-vector space of odd dimension $N\ge1$ with a given 
quadratic form such that the associated symplectic form has radical $\rr$ of dimension $1$ and the 
restriction of the quadratic form to $\rr$ is nonzero. In this case there is an obvious morphism from $G$ to
to the symplectic group $G'$ of $V/\rr$ which is an isomorphism of abstract groups. From the definitions we
see that this morphisms maps each stratum of $G$ bijectively onto a stratum of $G'$ (which has been described
in 3.5).

\subhead 3.9\endsubhead
For an even $N\in\NN$ let ${}'Z^2_N$ be the set of partitions with labels $\nu=(\nu_1\ge\nu_2\ge\do)$ in 
$Z^2_N$ (see 3.4) such that the number of nonzero entries of $\nu$ is even.

(a) {\it If $N$ is even, then there is a canonical bijection ${}'Z^2_N\lra BP^{N/2}_{0,4}$.}
\nl
To $\nu\in{}'Z^2_N$ we associate $\l=(\l_1,\l_2,\l_3,\do)$ as follows: each string $2a,2a,2a,\do$ of $\nu$ 
of odd length or of even length and label $1$ is replaced by $a-1,a+1,a-1,a+1,\do$ of the same length (if 
the string has even origin) or $a+1,a-1,a+1,a-1,\do$ of the same length (if the string has odd origin); each 
string $2a,2a,2a,\do$ of $\nu$ of even length and label $0$ is replaced by $a-2,a+2,a-2,a+2,\do$ of the same 
length (if the string has even origin) or $a,a,a,a,\do$ of the same length (if the string has odd origin); 
each string $2a+1,2a+1,2a+1,\do$ of $\nu$ (necessarily of even length) is replaced by $a-1,a+2,a-1,a+2,\do$ 
of the same length (if the string has even origin) or $a+1,a,a+1,a,\do$ of the same length (if the string 
has odd origin). The resulting entries form a bipartition $\l\in BP^{N/2}_{0,4}$. Now $\nu\m\l$ establishes 
the bijection (a).

Assume for example that $N=8$. The bijection (a) is:

$(62\do)\lra(40\do)$

$((44)_1\do)\lra(31\do)$ 

$((44)_0\do)\lra(22\do)$

$(4211\do)\lra(3010\do)$

$(3311\do)\lra(2110\do)$

$((2222)_1\do)\lra(2020\do)$

$((2222)_0\do)\lra(1111\do)$

$((22)_11111\do)\lra(201010\do)$

$((22)_01111\do)\lra(111010\do)$

$(11111111\do)\lra(10101010\do)$
\nl
Here we write $\do$ instead of $000\do$. (Compare \cite{\LSS, 6.2}.)

\subhead 3.10\endsubhead 
Assume that $G=SO(V)$ where $V$ is a $\kk$-vector space of even dimension $N$ with a fixed nondegenerate
quadratic form. Let $g\in G$. For any $x\in\kk^*$ let $V_x$ be the generalized $x$-eigenspace of $g:V@>>>V$.
Let $d_x=\dim V_x$. 
For any $x\in\kk^*$ such that $x^2\ne1$ let $\l^x_1\ge\l^x_2\ge\l_3^x\ge\do$ be the partition whose nonzero
terms are the sizes of the Jordan blocks of $x\i g:V_x@>>>V_x$. 
For $x\in\kk^*$ such that $x^2=1$ let $\nu^x\in {}'Z^1_{d_x}$ (if $p\ne2$) and
$\nu^x\in {}'Z^2_{d_x}$ (if $p=2$) be again the partition of $d_x$ 
whose nonzero terms are the sizes of the Jordan blocks of the unipotent element $x\i g\in SO(V_x)$. (When 
$p=2$, $\nu^x$ should also include a labelling with $0$ and $1$ associated to $x\i g$ viewed as an
element of $Sp(V_x)$ as in \cite{\WEUII, 1.4}.) Let $\l^x=(\l^x_1,\l^x_2,\l_3^x,\do)$
be the bipartition of $d_x/2$ associated to $\nu^x$ by 3.6(b), 3.9(a). Thus
$\l^x\in BP^{d_x/2}_{0,2}$ (if $p\ne2$), $\l^x\in BP^{d_x/2}_{0,4}$ (if $p=2$).
Note that $\l^x$ is the bipartition such that the Springer representation attached to the unipotent element 
$x\i g\in SO(V_x)$ (an irreducible representation of the Weyl group of type $D_{d_x/2}$) is indexed by 
$\l^x$. 
Define ${}^g\l=({}^g\l_1,{}^g\l_2,{}^g\l_3,\do)$ by ${}^g\l_j=\sum_x\l^x_j$ where $x$ runs over
a set of representatives for the orbits of the involution $a\m a\i$ of $\kk^*$. Note 
that ${}^g\l\in BP^{N/2}_{0,4}$ and that $g\m{}^g\l$ defines a (surjective) map $G@>>>BP^{N/2}_{0,4}$. From 
the definitions we see that the fibres of this map are exactly the strata of $G$ (except for the fibre over 
a bipartition $(\l_1,\l_2,\l_3,\do)$ with $\l_1=\l_2,\l_3=\l_4,\do$ in which case the fibre is a union of 
two strata). If $g\in G$ and ${}^g\l=(\l_1,\l_2,\l_3,\do)$ then 
$$\dim(\cb_g)=\sum_{k\ge1}((N/2)-(\l_1+\l_2+\do+\l_k)).\tag a$$ 

Viewing $W$ as a subgroup of index $2$ of a Weyl group $W'$ of type $B_n$, we can associate to any 
$\l\in BP^{N/2}$ one or two irreducible representations of $W$ which appear in the restriction to $W$ of the
irreducible representation of $W'$ indexed by $\l$; the representation(s) of $W$ associated to 
$\l=(\l_1,\l_2,\l_3,\l_4,\do)$ are the same as those associated to $\io(\l):=(\l_2,\l_1,\l_4,\l_3,\do)$; here
$\io:BP^{N/2}@>>>BP^{N/2}$ is an involution with set of orbits denoted by $BP^{N/2}/\io$. This gives a 
surjective map $f:\Irr(W)@>>>BP^{N/2}/\io$ whose fibre at the orbit of $\l$ has one element if 
$\l\ne\io(\l)$ and two elements if $\l=\io(\l)$. Let $\io':\Irr(W)@>>>\Irr(W)$ be the involution
whose orbits are the fibres of $f$ and let $\cs_2(W)/\io'$ be the set of orbits of the restriction of $\io'$
to $\cs_2(W)$. The results in this subsection show that $f$ induces a bijection
$$\cs_2(W)/\io'@>\si>>BP^{N/2}_{0,4}.\tag b$$
We have used the fact that the intersection of $BP^{N/2}_{0,4}$ with an orbit of $\io:BP^{N/2}@>>>BP^{N/2}$ 
has at most one element; more precisely,
$$\{\l\in BP^{N/2};\l\in BP^{N/2}_{0,4}\text{ and } \io(\l)\in BP^{N/2}_{0,4}\}
=\{\l\in BP^{N/2};\l=\io(\l)\}.$$
Under the identification (b), the map $g\m{}^g\l$, $G@>>>BP^{N/2}_{0,4}$ becomes the map $g\m E$ (up to the 
action of $\io'$) where $g\in G_E$.

\subhead 3.11\endsubhead 
Assume that $p\ne2$ and $n\ge3$. If $G=SO_{2n+1,\kk}$ then the stratum of minimal dimension $>0$ consists of
a semisimple class of dimension $2n$; if $G=Sp_{2n,\kk}/\pm1$ then the stratum of minimal dimension $>0$ 
consists of a unipotent class of dimension $2n$ (that of transvections). The corresponding $E\in\Irr(W)$ is
one dimensional.

\subhead 3.12\endsubhead 
Assume that $G$ is simple of type $E_8$. In this case $G$ has exactly $75$ strata. If $p\ne2,3$ then exactly
$70$ strata contain unipotent elements. If $p=2$ (resp. $p=3$) then exactly $74$ (resp. $71$) strata contain
unipotent elements. The unipotent class of dimension $58$ is a stratum. If $p\ne2$, there is a stratum 
which is a union of a semisimple class and a unipotent class (both of dimension $128$); in particular this 
stratum is disconnected. 

\head 4. A map from conjugacy classes in $W$ to $2$-special representations of $W$\endhead
\subhead 4.1\endsubhead 
In this subsection we shall define a canonical surjective map
$${}'\Ph:cl(W)@>>>\cs_2(W).\tag a$$
We preserve the setup of 2.5. We will first define the map (a) assuming that $G$ is simple.
In \cite{\WEU} we have defined for any $r\in\cp$ a surjective map $cl(W)@>>>\cu^r$; we denote this map by 
$\Ph^r$. Let $C\in cl(W)$. We define an element $\Ph(C)\in\cu^*$ as follows. 
If $\Ph^r(C)\in h^r(z_r)$ (with $z_r\in\cu^0$) for all $r\in\cp$ then $z_r=z$ is independent of $r$ (see 
\cite{\WEUII, 0.4}) and we define $\Ph(C)$ to be the equivalence class of $h^r(z)$ for any $r\in\cp$. If 
$\Ph^r(C)\n h^r(\cu^0)$ for some $r\in\cp$ then $r$ is unique. (The only case where $r$ can be possibly not 
unique is in type $E_8$ in which case we use the tables in \cite{\WEUII, 2.6}.) We then define $\Ph(C)$ to 
be the equivalence class of $\Ph^r(C)$. Thus we have defined a surjective map $\Ph:cl(W)@>>>\cu^*$. 
By composing $\Ph^r$ with $\ps^r:\cu^r@>\si>>\cs^r_2(W)$, see 2.5, and with the inclusion
$\cs^r_2(W)\sub\cs_2(W)$, we obtain a map ${}'\Ph^r:cl(W)@>>>\cs_2(W)$. Similarly, by composing $\Ph$ with 
$\ps^*:\cu^*@>\si>>\cs_2(W)$, see 2.5(a), we obtain a surjective map ${}'\Ph:cl(W)@>>>\cs_2(W)$.
Note that for $C\in cl(W)$, ${}'\Ph(C)$ can be described as follows. If ${}'\Ph^r(C)\in\cs_1(W)$ for all 
$r\in\cp$ then ${}'\Ph^r(C)$ is indepedent of $r$ and we have ${}'\Ph(C)={}'\Ph^r(C)$ for any $r$. If 
${}'\Ph^r(C)\n\cs_1(W)$ for some $r\in\cp$ then such $r$ is unique and we have ${}'\Ph(C)={}'\Ph^r(C)$.

We return to the general case. We write the adjoint group of $G$ as a product $\prod_iG_i$ where each $G_i$ 
is simple with Weyl group $W_i$. We can identify $W=\prod_iW_i$, $cl(W)=\prod_icl(W_i)$, 
$\cs_2(W)=\prod_u\cs_2(W_i)$ (via external tensor product). Then ${}'\Ph_i:cl(W_i)@>>>\cs_2(W_i)$ is defined
as above for each $i$. We set ${}'\Ph=\prod_i{}'\Ph_i:cl(W)@>>>\cs_2(W)$.

For $C,C'$ in $cl(W)$ we write $C\si C'$ if ${}'\Ph(C)={}'\Ph(C')$. This is an equivalence relation on
$cl(W)$. Let $\ucl(W)$ be the set of equivalence classes. Note that:

(b) {\it ${}'\Ph$ induces a bijection $\ucl(W)@>>>\cs_2(W)$.}
\nl
We see that, via (b), 

(c) {\it the strata of $G$ are naturally indexed by the set $\ucl(W)$.}

\subhead 4.2\endsubhead 
We preserve the setup of 2.5. Now ${}'\Ph$ in 4.1(a) is a map between two sets which depend only on $W$, not
on the underlying root system, see 1.1(b). We show that 

(a) {\it ${}'\Ph$ itself depends only on $W$, not on the underlying root system.}
\nl
We can assume that $G$ is adjoint, simple. We can also assume that $G$ is not of simply laced type. In this 
case there is a unique $r\in\cp$ such that $\cs_2(W)=\cs^r_2(W)$ so that we have simply
${}'\Ph={}'\Ph^r:cl(W)@>>>\cs_2(W)$. Thus ${}'\Ph$ is the composition 

(b) $cl(W)@>\Ph^r>>\cu^r@>\ps^r>>\cs_2(W)$.
\nl
We now use the fact the maps in (b) are compatible with the exceptional isogeny between groups $G^2$ of type
$B_n$ and $C_n$ or of type $F_4$ and $F_4$ (resp. between groups $G^3$ of type $G_2$ and $G_2$). This 
implies (a).

\subhead 4.3\endsubhead 
Assume that $G$ is simple. The map ${}'\Ph$ in 4.1 is defined in terms of ${}'\Ph^r$ which is the 
composition of $\Ph^r:cl(W)@>>>\cu^r$ (which is described explicitly in each case in \cite{\WEUII}) and
$\ps^r:\cu^r@<<<\cs^r_2(W)$ which is given by the Springer correspondence. Therefore ${}'\Ph$ is explicitly
computable. In this subsection we describe this map in the case where $W$ is of classical type.

If $W$ is of type $A_n$, $n\ge1$, then $cl(W)$ can be identified with the set of partitions of $n$: to a
conjugacy class of a permutation of $n$ objects we associate the partition whose nonzero terms are the sizes
of the disjoint cycles of which the permutation is a product. We identify $\cs_2(W)=\Irr(W)$ with the set of
partitions in the standard way (the unit representation corresponds to the partition $(n,0,0\do)$). With 
these identifications the map ${}'\Ph$ is the identity map.

Assume now that $W$ is a Weyl group of type $B_n$ or $C_n$, $n\ge2$. Let $X$ be a set with $2n$ elements with
a given fixed point free involution $\t$. We identify $W$ with the group of permutations of $X$ which 
commute with $\t$. To any $w\in W$ we can associate an element $\nu\in Z^2_{2n}$ (see 3.4) as follows. The
nonzero terms of the partition $\nu$ are the sizes of the disjoint cycles of which $w$ is a product. To each
string $c,c,\do,c$ of $\nu$ of even length with $c>0$ even we attach the label $1$ if at least one of its 
terms represents a cycle which commutes with $\t$; otherwise we attach to it the label $0$. This defines a
(surjective) map $cl(W)@>>>Z^2_{2n}$ which by results of \cite{\WEUII} can be identified with the map
$\Ph^2:cl(W)@>>>\cu^2$. Composing this with the bijection 3.4(b) we obtain a surjective map
$cl(W)@>>>BP^n_{2,2}$ or equivalently (see 3.5(b)) $cl(W)@>>>\cs_2(W)$. This is the same as ${}'\Ph$.

Next we assume that $W$ is a Weyl group of type $D_n$, $n\ge4$. We can identify $W$ with the 
group of {\it even} permutations of $X$ (as above) which commute with $\t$ (as above).
To any $w\in W$ we associate an element $\nu\in Z^2_{2n}$ as for type $B_n$ above. This element is actually
contained in ${}'Z^2_{2n}$ (see 3.9) since $w$ is an even permutation. This defines a
(surjective) map $cl(W)@>>>{}'Z^2_{2n}$ which by results of \cite{\WEUII} can be identified with the 
composition of $\Ph^2:cl(W)@>>>\cu^2$ with the obvious map from $\cu^2$ to the set of orbits of the 
conjugation action of the full orthogonal group on $\cu^2$. Composing this with the bijection 3.9(a) we 
obtain a surjective map $cl(W)@>>>BP^n_{0,4}$ or equivalently (see 3.10(b)) a surjective map 
$cl(W)@>>>\cs_2(W)/\io'$ (notation of 3.10). This is the same as the composition of ${}'\Ph$ with the 
obvious map $\cs_2(W)@>>>\cs_2(W)/\io'$.

\subhead 4.4\endsubhead 
In this and the next four subsections we describe the map ${}'\Ph:cl(W)@>>>\cs_2(W)$ in the case where $W$ 
is of exceptional type. The results will be expressed as diagrams $[a,b,\do]\m d_n$ where $a,b,\do$ is the 
list of conjugacy classes in $W$ (with notation of \cite{\CA}) which are mapped by ${}'\Ph$ to an irreducible
representation $E$ denoted $d_n$ (here $d$ denotes the degree of $E$ and the index $n=n_E$ as in 0.2). We
also mark by $*_r$ those $E$ which are in $\cs_2(W)-\cs_1(W)$; here $r$ is the unique prime such that
$E\in\cs^r_2(W)$. Note that the notation $d_n$ does not determine
$E$ for types $G_2$ and $F_4$; for these types it may happen that there are two $E$'s with same $d_n$.

\mpb

{\it Type $G_2$.}

$[G_2]\m1_0$ 

$[A_2]\m2_1$ 

$[A_1+\tA_1]\m2_2$

$[\tA_1]\m1_3, *_3$

$[A_1]\m1_3$ 

$[A_0]\m1_6$ 

\subhead 4.5\endsubhead 
{\it Type $F_4$.}

$[F_4]\m1_0$

$[B_4]\m4_1$     

$[F_4(a_1)]\m9_2$      

$[D_4,B_3]\m8_3$       

$[C_3+A_1,C_3]\m8_3$         

$[D_4(a_1)]\m12_4$       

$[A_3+\tA_1]\m16_5$  

$[A_3]\m9_6$       

$[B_2+A_1]\m9_6, *_2$ 

$[\tA_2+\tA_2]\m6_6$    

$[A_2+\tA_1]\m4_7$    

$[\tA_2+A_1]\m 4_7,*_2$

$[B_2]\m 4_8,*_2$     

$[\tA_2]\m 8_9$        

$[A_2]\m 8_9$ 

$[4A_1,3A_1,2A_1+\tA_1,A_1+\tA_1]\m 9_{10}$     

$[2A_1]\m 4_{13}$         

$[A_1]\m 2_{16}$       

$[\tA_1]\m 2_{16},*_2$

$[A_0]\m 1_{24}$           

\subhead 4.6\endsubhead              
{\it Type $E_6$.}

$[E_6]\m1_0$

$[E_6(a_1)]\m6_1$

$[D_5]\m20_2$

$[E_6(a_2)]\m 30_3$    

$[A_5+A_1,A_5]\m15_4$ 

$[D_5(a_1)]\m64_4$

$[A_4+A_1]\m60_5$

$[D_4]\m 24_6 $

$[A_4]\m81_6$

$[D_4(a_1)]\m 80_7$

$[A_3+2A_1,A_3+A_1]\m60_8$

$[3A_2,2A_2+A_1]\m 10_9$

$[A_3]\m 81_{10}$

$[A_2+2A_1]\m 60_{11}$

$[2A_2]\m 24_{12}$

$[A_2+A_1]\m 64_{13}$

$[A_2]\m 30_{15}$

$[4A_1,3A_1]\m 15_{16}$     

$[2A_1]\m 20_{20}$      

$[A_1]\m 6_{25}$

$[A_0]\m 1_{36}$     

\subhead 4.7\endsubhead
{\it Type $E_7$.}

$[E_7]\m 1_0$

$[E_7(a_1)]\m 7_1$

$[E_7(a_2)]\m 27_2$

$[E_7(a_3)]\m 56_3$

$[E_6]\m 21_3$

$[E_6(a_1)]\m120_4$

$[D_6+A_1,D_6]\m35_4$

$[A_7]\m189_5$

$[A_6]\m105_6$

$[D_6(a_1)]\m210_6$

$[D_5+A_1]\m168_6$            

$[E_7(a_4)]\m315_7$

$[D_5]\m189_7$

$[E_6(a_2)]\m405_8$

$[D_6(a_2)+A_1,D_6(a_2)]\m280_8$

$[A_5+A_2,(A_5+A_1)']\m70_9$

$[(A_5+A_1)'',A''_5]\m216_9$

$[D_5(a_1)+A_1]\m378_9$

$[D_5(a_1)]\m420_{10}$

$[A_4+A_2]\m210_{10}$

$[A_4+A_1]\m512_{11}$

$[A'_5]\m105_{12}$

$[D_4+3A_1,D_4+2A_1,D_4+A_1]\m84_{12}$

$[A_4]\m420_{13}$

$[2A_3+A_1,A_3+A_2+A_1]\m210_{13}$

$[A_3+A_2]\m378_{14}$       

$[D_4]\m105_{15}$             

$[D_4(a_1)+A_1]\m405_{15}$           

$[A_3+A_2]\m84_{15},*_2$   

$[A_3+3A_1,(A_3+2A_1)']\m216_{16}$      

$[D_4(a_1)]\m315_{16}$                  

$[(A_3+2A_1)'',(A_3+A_1)'']\m280_{17}$

$[3A_2,2A_2+A_1]\m70_{18}$

$[(A_3+A_1)']\m189_{20}$

$[A_3]\m210_{21}$

$[2A_2]\m168_{21}$

$[A_2+3A_1]\m105_{21}$

$[A_2+2A_1]\m189_{22}$

$[A_2+A_1]\m120_{25}$

$[7A_1,6A_1,5A_1,(4A_1)']\m15_{28}$   

$[A_2]\m56_{30}$

$[(4A_1)'',(3A_1)'']\m35_{31}$

$[(3A_1)']\m21_{36}$

$[2A_1]\m27_{37}$

$[A_1]\m7_{46}$

$[A_0]\m1_{63}$

\subhead 4.8\endsubhead
{\it Type $E_8$.}

$[E_8]\m1_0$

$[E_8(a_1)]\m8_1$

$[E_8(a_2)]\m35_2$

$[E_8(a_4)]\m112_3$

$[E_7+A_1,E_7]\m84_4$

$[E_8(a_5)]\m210_4$ 

$[D_8]\m560_5$

$[E_7(a_1)]\m567_6$

$[E_8(a_3)]\m700_6$

$[D_8(a_1),D_7]\m400_7$

$[E_8(a_7)]\m1400_7$

$[E_8(a_6)]\m1400_8$

$[E_7(a_2)+A_1,E_7(a_2)]\m1344_8$

$[E_6+A_2,E_6+A_1]\m448_9$

$[D_8(a_2)]\m3240_9$

$[D_7(a_1)]\m1050_{10},*_2$

$[A''_7]\m175_{12},*_3$    

$[A_8]\m2240_{10}$

$[E_7(a_3)]\m2268_{10}$ 

$[E_6(a_1)+A_1]\m4096_{11}$

$[D_8(a_3)]\m1400_{11}$

$[E_6]\m525_{12}$

$[D_7(a_2)]\m4200_{12}$

$[D_6+2A_1,D_6+A_1,D_6]\m972_{12}$

$[E_6(a_1)]\m2800_{13}$

$[A_7+A_1]\m4536_{13}$

$[A'_7]\m6075_{14}$

$[A_6+A_1]\m2835_{14}$

$[D_5+A_2]\m840_{14},*_2$  

$[A_6]\m4200_{15}$

$[D_6(a_1)]\m5600_{15}$

$[E_8(a_8)]\m4480_{16}$

$[D_5+2A_1,D_5+A_1]\m3200_{16}$

$[E_7(a_4)+A_1,E_7(a_4)]\m7168_{17}$

$[2D_4,D_6(a_2)+A_1,D_6(a_2)]\m4200_{18}$

$[E_6(a_2)+A_2,E_6(a_2)+A_1]\m3150_{18}$

$[A_5+A_2+A_1,A_5+A_2,A_5+2A_1,(A_5+A_1)'']\m2016_{19}$

$[D_5(a_1)+A_3,D_5(a_1)+A_2]\m1344_{19}$

$[D_5]\m2100_{20}$

$[2A_4,A_4+A_3]\m420_{20}$

$[E_6(a_2)]\m5600_{21}$

$[D_4+A_3]\m4200_{21}$

$[(A_5+A_1)']\m3200_{22}$

$[D_5(a_1)+A_1]\m6075_{22}$

$[A_4+A_2+A_1]\m2835_{22}$

$[A_4+A_2]\m4536_{23}$

$[A_4+2A_1]\m4200_{24}$

$[D_4+A_2]\m168_{24},*_2$        

$[D_5(a_1)]\m2800_{25}$

$[A_4+A_1]\m4096_{26}$

$[2D_4(a_1),D_4(a_1)+A_3,(2A_3)'']\m840_{26}$

$[D_4+4A_1,D_4+3A_1,D_4+2A_1,D_4+A_1]\m700_{28}$

$[D_4(a_1)+A_2]\m2240_{28}$

$[2A_3+2A_1,A_3+A_2+2A_1,2A_3+A_1,A_3+A_2+A_1]\m1400_{29}$

$[A_4]\m2268_{30}$

$[(2A_3)']\m3240_{31}$

$[D_4(a_1)+A_1]\m1400_{32}$

$[A_3+A_2]\m972_{32},*_2$      

$[A_3+4A_1,A_3+3A_1,(A_3+2A_1)'']\m1050_{34}$

$[D_4]\m525_{36}$

$[4A_2,3A_2+A_1,2A_2+2A_1]\m175_{36}$

$[D_4(a_1)]\m1400_{37}$

$[(A_3+2A_1)',A_3+A_1]\m1344_{38}$

$[3A_2,2A_2+A_1]\m448_{39}$                       

$[2A_2]\m700_{42}$

$[A_2+4A_1,A_2+3A_1]\m400_{43}$

$[A_3]\m567_{46}$

$[A_2+2A_1]\m560_{47}$

$[A_2+A_1]\m210_{52}$

$[8A_1,7A_1,6A_1,5A_1,(4A_1)'']\m50_{56}$ 

$[A_2]\m112_{63}$

$[(4A_1)',3A_1]\m84_{64}$

$[2A_1]\m35_{74}$

$[A_1]\m8_{91}$

$[A_0]\m1_{120}$                      

\subhead 4.9\endsubhead
In the tables in 4.4-4.8 the $E$ which are not marked with $*_r$ are in $\cs_1(W)$; they are expressed
explicitly in the form $j_{W_{e'}}^W(E')$ with $e'\in V^*$, $E'\in\cs(W_{e'})$ in the tables of 
\cite{\SPEC}. 

We now consider the $E$ in the tables 4.4-4.8 which are marked with $*_r$.

\mpb

Type $G_2$: 

$1_3=j_{W'}^W(\text{sign})$ where $W'$ is of type $A_2$ but not of form $W_{e'},e'\in V^*$.

\mpb

Type $F_4$: 

$9_6=j_{W'}^W(\text{E'})$ where $W'$ is of type $B_4$ but not of form $W_{e'},e'\in V^*$ and $\dim E'=6$,
$n_{E'}=6$;

$4_7=j_{W'}^W(\text{sign})$ where $W'$ is of type $A_3A_1$ but not of form $W_{e'},e'\in V^*$;

$4_8=j_{W'}^W(\text{sign})$ where $W'$ is of type $B_2B_2$;

$2_{12}=j_{W'}^W(\text{sign})$ where $W'$ is of type $B_4$ but not of form $W_{e'},e'\in V^*$.

\mpb

Type $E_7$:

$84_{15}=j_{W'}^W(\text{sign})$ where $W'$ is of type $D_4A_1A_1A_1$.

\mpb

Type $E_8$:

$1050_{10}=j_{W'}^W(\text{E'})$ where $W'$ is of type $D_6A_1A_1$ and $\dim E'=30$, $n_{E'}=10$;

$175_{12}=j_{W'}^W(\text{sign})$ where $W'$ is of type $A_2A_2A_2A_2$.

$840_{14}=j_{W'}^W(\text{sign})$ where $W'$ is of type $A_3A_3A_1A_1$.

$168_{24}=j_{W'}^W(\text{sign})$ where $W'$ is of type $D_4D_4$.

$972_{32}=j_{W'}^W(\text{sign})$ where $W'$ is of type $D_6A_1A_1$.

\subhead 4.10\endsubhead
For any $C\in cl(W)$ let $m_C$ be the dimension of the $1$-eigenspace
of an element in $C$ in the reflection representation of $W$. We have the following result.

(a) {\it For any $E\in\cs_2(W)$, the restriction of $C\m m_C$ to ${}'\Ph\i(E)\sub cl(W)$ reaches its 
minimum at a unique element of ${}'\Ph\i(E)$, denoted by $C_E$.}
\nl
We can assume that $G$ is simple. When $G$ is of exceptional type, (a) follows from the tables 4.4-4.8.
When $G$ is of classical type, (a) follows from \cite{\WEUII,0.2}.

Note that $E\m C_E$ is a cross section of the surjective map ${}'\Ph:cl(W)@>>>\cs_2(W)$. It defines a
bijection of $\cs_2(W)$ with a subset $cl_0(W)$ of $cl(W)$.

\head 5. A second approach\endhead
\subhead 5.1\endsubhead
In this section we sketch another approach to defining the strata of $G$ in which Springer representations 
do not appear. Let $cl(G)$ be the set of conjugacy classes in $G$.
Let $\ul:W@>>>\NN$ be the length function of the Coxeter group $W$. For $w\in W$ let 
$$G_w=\{g\in G;(B,gBg\i)\in\co_w\text{ for some }B\in\cb\}.$$
For $C\in cl(W)$ let 
$$C_{min}=\{w\in C;\ul:C@>>>\NN\text{ reaches minimum at }w\}$$
and let $G_C=G_w$ where $w\in C_{min}$.

As pointed out in \cite{\WEU, 0.2}, from \cite{\WEU, 1.2(a)} and \cite{\GP, 8.2.6(b)} it follows that $G_C$ 
is independent of the choice of $w$ in $C_{min}$. From \cite{\WEU} it is known that $G_C$ contains unipotent
elements; in particular, $G_C\ne\emp$. Clearly, $G_C$ is a union of conjugacy classes. Let 
$$\d_C=\min_{\g\in cl(G);\g\sub G_C}\dim\g,$$
$$\bx{G_C}=\cup_{\g\in cl(G);\g\sub G_C,\dim\g=\d_C}\g.$$
Then $\bx{G_C}$ is $\ne\emp$, a union of conjugacy classes of fixed dimension, $\d_C$. We have the following
result.

\proclaim{Theorem 5.2}Let $C\in cl(W)$, $E\in\cs_2(W)$ be such that ${}'\Ph(C)=E$, see 4.1. We have 
$\bx{G_C}=G_E$.
\endproclaim
We can assume that $G$ is almost simple and that $\kk$ is an algebraic closure of a finite field. The proof 
in the case of exceptional groups is reduced in 5.3 to a computer calculation. The proof for classical 
groups, which is based on combining the techniques of \cite{\WEU}, \cite{\CSM} and \cite{\DIST}, will be 
given elsewhere. 

\subhead 5.3\endsubhead
In this subsection we assume that $\kk$ is an algebraic closure of a finite field $\FF_q$ and that $G$ is 
simply connected, defined and split over $\FF_q$ with Frobenius map $F:G@>>>G$. Let $\g$ be an $F$-stable 
conjugacy class of $G$. Let $\g'=\{g_s;g\in\g\}$, an $F$-stable semisimple conjugacy class in $G$. For every
$s\in\g'$ let $\g(s)=\{u\in Z_G(s);u\text{ unipotent, }us\in\g\}$, a unipotent conjugacy class of $Z_G(s)$. 
We fix $s_0\in\g'{}^F$ and we set $H=Z_G(s_0)$, $\g_0=\g(s_0)$. Let $W_H$ be the Weyl group of $H$. As in 
2.1, we can regard $W_H$ as a subgroup of $W$ (the imbedding of $W_H$ into $W$ is canonical up to composition
with an inner automorphism of $W$). By replacing if necessary $F$ by a power of $F$, we can assume that $H$ 
contains a maximal torus which is defined and split over $\FF_q$. For any $F$-stable maximal torus $T$ of 
$G$, $R_T^1$ is the virtual representation of $G^F$ defined as in \cite{\DL, 1.20} (with $\th=1$ and with 
$B$ omitted from notation). Replacing $T,G$ by $T',H$ where $T'$ is an $F$-stable maximal torus of $H$, we 
obtain a virtual representation $R_{T',H}^1$ of $H^F$. 
For any $z\in W$ we denote by $R_z^1$ the virtual representation $R_T^1$ of $G^F$ where $T$ is an $F$-stable
maximal torus of $G$ of type given by the conjugacy class of $z$ in $W$. For any $z'\in W_H$ we denote by 
$R_{z',H}^1$ the virtual representation $R_{T',H}^1$ of $H^F$ where $T'$ is an $F$-stable maximal torus of 
$H$ of type given by the conjugacy class of $z'$ in $W_H$. For $E'\in\Irr W$ we set 
$R_{E'}=|W|\i\sum_{y\in W}\tr(y,E')R_y^1$. Then for any $z\in W$ we have 
$R_z^1=\sum_{E'\in\Irr W}\tr(z,E')R_{E'}$.

Let $w\in W$. We show:
$$\align&|\{(g,B)\in\g^F\T\cb^F;(B,gBg\i)\in\co_w\}|\\&=
|G^F||H^F|\i\sum_{E\in\Irr W,E'\in\Irr W,E''\in\Irr W_H,y}\tr(T_w,E_q)(\r_E,R_{E'})\\&
\T(E'|_{W_H}:E'')|Z_{W_H}(y)|\i tr(y,E'')\sum_{u\in\g_0^F}\tr(u,R_{y,H}^1)\tag a\endalign$$
where $y$ runs over a set of representatives for the conjugacy classes in $W_H$ and $T_w,E_q,\r_E$ are as in
\cite{\WEU, 1.2}. Let $N$ be the left hand side of (a). As in \cite{\WEU, 1.2(c)} we see that
$$N=\sum_{E\in\Irr W}\tr(T_w,E_q)A_E$$
with
$$A_E=|G^F|\i\sum_{g\in\g^F}\sum_T|T^F|(\r_E,R_T^1)\tr(g,R_T^1)$$
where $T$ runs over all maximal tori of $G$ defined over $\FF_q$. We have
$$\align&A_E=|G^F|\i\sum_{s\in\g'{}^F,u\in\g(s)^F}\sum_T|T^F|(\r_E,R_T^1)\tr(su,R_T^1)\\&
=|H^F|\i\sum_{u\in\g_0^F}\sum_T|T^F|(\r_E,R_T^1)\tr(s_0u,R_T^1).\endalign$$
By \cite{\DL, 4.2} we have
$$\tr(s_0u,R_T^1)=|H^F|\i\sum_{x\in G^F;x\i Tx\sub H}\tr(u,R_{x\i Tx,H}^1)$$
hence
$$\align&A_E=|H^F|^{-2}\sum_{u\in\g_0^F}\sum_T|T^F|(\r_E,R_T^1)
\sum_{x\in G^F;x\i Tx\sub H}\tr(u,R_{x\i Tx,H}^1)\\&
=|G^F||H^F|^{-2}\sum_{T'\sub H}|T'{}^F|(\r_E,R_{T'}^1)\sum_{u\in\g_0^F}\tr(u,R_{T',H}^1)\endalign$$
where $T'$ runs over the maximal tori of $H$ defined over $\FF_q$. Using the classification of maximal tori
of $H$ defined over $\FF_q$ we obtain
$$\align&A_E=|G^F||H^F|\i|W_H|\i\sum_{z\in W_H}(\r_E,R_z^1)\sum_{u\in\g_0^F}\tr(u,R_{z,H}^1)\\&
=|G^F||H^F|\i|W_H|\i\sum_{z\in W_H}\sum_{E'\in\Irr W}\tr(z,E')(\r_E,R_{E'})\sum_{u\in\g_0^F}\tr(u,R_{z,H}^1).
\endalign$$
This clearly implies (a).

Now assume that $G$ is almost simple of exceptional type and that $w$ has minimal length in its conjugacy 
class in $W$. We can also assume that $q-1$ is sufficiently divisible. Then the right hand side of (a) can 
be explicitly determined using a computer. Indeed, it is an entry of the product of several large matrices 
whose entries are explicitly known. In particular the quantities $\tr(T_w,E_q)$ (known from the works of 
Geck and Geck-Michel, see \cite{\GP, 11.5.11}) are available through the CHEVIE package \cite{\GH}. The 
quantities $(\r_E,R_{E'})$ are coefficients of the nonabelian Fourier transform in \cite{\ORA, 4.15}. The 
quantities $(E'|_{W_H}:E'')$ are available from the induction tables in the CHEVIE package. The quantities 
$\tr(y,E'')$ are available through the CHEVIE package. The quantities $\tr(u,R_{y,H}^1)$ are Green 
functions; I thank Frank L\"ubeck for providing to me tables of Green functions for groups of rank $\le8$ in
GAP format. I also thank Gongqin Li for her help with programming in GAP to perform the actual computation 
using these data.

Thus the number $|\{(g,B)\in\g^F\T\cb^F;(B,gBg\i)\in\co_w\}|$ is explicitly computable. It turns out that it
is a polynomial in $q$. Note that $\{(g,B)\in\g\T\cb;(B,gBg\i)\in\co_w\}$ is nonempty if and only if this 
polynomial is non zero. Thus the condition that $\g\sub G_w$ can be tested. This can be used to check that 
Theorem 5.2 holds for exceptional groups.

\subhead 5.4\endsubhead
If $C$ is the conjugacy class containing the Coxeter elements of $W$ then $G_C=\bx{G_C}$ is the union of all
conjugacy classes of dimension $\dim G-\rk(G)$, see \cite{\ST}. 

\head 6. Variants\endhead
\subhead 6.1\endsubhead
The results in this subsection will be proved elsewhere.
In this subsection we assume that $G$ is simple and that $G'$ is a disconnected reductive algebraic group 
$G$ over $\kk$ with identity component $G$ such that $G'/G$ is cyclic of order $r$ and
such that the homomorphism $\e:G'/G@>>>\Aut(W)$ (the automorphism group of $W$ as a Coxeter group)
induced by the conjugation action of $G'/G$ on $G$ is injective.
Note that $(G,r)$ must be of type $(A_n,2)$ ($n\ge2$) or $(D_n,2)$ ($n\ge4$) or $(D_4,3)$ or $(E_6,2)$. Let 
$D$ be a connected component of $G'$ other than $G$. We will give a 
definition of the strata of $D$, extending the definition of strata of $G$. 
Let $\e_D:W@>>>W$ be the image of $D$ under $\e$.
Let $cl_DW$ be the set of conjugacy classes in $W$ twisted by $\e_D$ (as in \cite{\DIST,0.1}).
Let $cl(D)$ be the set of $G$-conjugacy classes in $D$. For $w\in W$ let 
$$D_w=\{g\in D;(B,gBg\i)\in\co_w\text{ for some }B\in\cb\}.$$
For $C\in cl_D(W)$ let 
$$C_{min}=\{w\in C;\ul:C@>>>\NN\text{ reaches minimum at }w\}.$$
and let $D_C=D_w$ where $w\in C_{min}$. This is independent of the choice of $w$ in $C_{min}$. One
can show that $D_C\ne\emp$. Clearly, $D_C$ is a union of $G$-conjugacy classes in $D$. Let 
$$\d_C=\min_{\g\in cl(D);\g\sub D_C}\dim\g,$$
$$\bx{D_C}=\cup_{\g\in cl(D);\g\sub D_C,\dim\g=\d_C}\g.$$
Then $\bx{D_C}$ is $\ne\emp$, a union of $G$-conjugacy classes of fixed dimension, $\d_C$. 
One can show that $\cup_{C\in cl_D(W)}\bx{D_C}=D$; moreover, one can show that if $C,C'\in cl_D(W)$, then
$\bx{D_C},\bx{D_{C'}}$ are either equal or disjoint. (Some partial results in this direction  are contained 
in \cite{\DIST}.) Let $\si$ be the equivalence relation on
$cl_D(W)$ given by $C\si C'$ if $\bx{D_C}=\bx{D_{C'}}$ and let $\ucl_D(W)$ be the set of equivalence 
classes. We see that there is a unique partition of $D$ into pieces (called {\it strata}) indexed by 
$\ucl_D(W)$ such that each stratum is of the form $\bx{D_C}$ for some $C\in cl_D(W)$.
One can show that the equivalence relation $\si$ on $cl_D(W)$ and the function $C\m d_C$ on $cl_D(W)$
depend only on $W$ and its automorphism $\e_D$; in particular they do not depend on $\kk$.
When $p=r$, each stratum of $D$ contains a unique unipotent $G$-conjugacy class in $D$; this gives a
bijection $\ucl_D(W)\lra\cu_D^r$ where $\cu_D^r$ is the set of unipotent $G$-conjugacy classes in $D$ 
(with $p=r$). This bijection coincides with the bijection $\ucl_D(W)\lra\cu_D^r$ described explicitly 
in \cite{\WEUIII}. Thus the strata of $D$ can also
be indexed by $\cu_D^r$. We can also index them by a certain set of irreducible representations of
$W^{\e_D}$ (the fixed point set of $\e_D:W@>>>W$) using the bijection \cite{\CDG, II} between $\cu_D^r$
and a set of irreducible representations of $W^{\e_D}$ (an extension of the Springer correspondence).

\subhead 6.2 \endsubhead
Assume that $G$ is adjoint.
We identify $\cb$ with the variety of Borel subalgebras of $\fg$. For any $\x\in\fg$ let
$\cb_\x=\{\fb\in\cb;\x\in\fb\}$ and let $d=\dim\cb_\x$. The subspace of $H_{2d}(\cb)$ spanned by the images 
of the fundamental classes of the irreducible components of $\cb_\x$ is an irreducible $W$-module denoted by
$\t_\x$. We also denote by $\t_\x$ the corresponding $W$-module over $\QQ$. Thus we have a well defined map 
$\fg@>>>\Irr W$, $\x\m\t_\x$. The nonempty fibres of this map are called the {\it strata} of $\fg$. Each 
stratum of $\fg$ is a union of adjoint orbits of fixed dimension; exactly one of these orbits is
nilpotent. The image of the map $\x\m\t_\x$ is the subset of $\Irr(W)$ denoted by $\ct^p_W$ in \cite{\REM};
when $p=0$ this is $\cs_1(W)$.

\subhead 6.3\endsubhead
In this subsection we assume that $G$ is semisimple simply connected. Let $K$ be the field of formal power
series $\kk((\e))$ and let $\hat G=G(K)$. Let $\hat\cb$ be the set of Iwahori subgroups of $\hat G$ viewed
as an increasing union of projective algebraic varieties over $\kk$. Let $\hat W$ be the affine Weyl group
associated to $\hat G$ viewed as an infinite Coxeter group. Let $G(K)_{rsc}$ be the set of all $g\in G(K)$
which are compact (that is such that $\hat\cb_g=\{B\in\hat\cb;g\in B\}$ is nonempty) and regular semisimple.
If $g\in G(K)_{rsc}$ then $\hat\cb_g$ a union of projective algebraic varieties of fixed dimension $d=d_g$
(see \cite{\KL} for a closely related result) hence the homology space $H_{2d}(\hat\cb_g)$ is well 
defined and it carries a natural $\hat W$-action (see \cite{\UNAC}). Similarly the homology space 
$H_{2d}(\hat\cb)$ is well defined and it carries a natural $\hat W$-action. The imbedding 
$h_g:\hat\cb_g@>>>\hat\cb$ induces a linear map $h_{g*}:H_{2d}(\hat\cb_g)@>>>H_{2d}(\hat\cb)$
which is compatible with the $\hat W$-actions. Hence $\hat W$ acts naturally on the (finite dimensional) 
subspace $E_g:=h_{g*}(H_{2d}(\hat\cb_g))$ of $H_{2d}(\hat\cb)$, but this action is not irreducible
in general. Note that $E_g$ is the subspace of $H_{2d}(\hat\cb)$ spanned by the images of the 
fundamental classes of the irreducible components of $\hat\cb_g,\bbq$ (we ignore Tate twists) hence is 
$\ne0$. For $g,g'\in G(K)_{rsc}$ we say that $g\si g'$ if $d_g=d_{g'}$ and $E_g=E_{g'}$. This is an 
equivalence relation on $G(K)_{rsc}$. The equivalence classes for $\si$ are called the {\it strata} of 
$G(K)_{rsc}$. Note that $G(K)_{rsc}$ is a union of countably many strata and each stratum is a union of 
conjugacy classes of $G(K)$ contained in $G(K)_{rsc}$.

\subhead 6.4\endsubhead
In this subsection we state a conjectural definition of the strata of $G$ in the case where $\kk=\CC$ based 
on an extension of a construction in \cite{\KL}. Let $K$ be as in 6.3. Let $g\in G$. Let $\fz\sub\fg$ be the
Lie algebra of $Z_G(g_s)$ and let $\x=\log(g_u)\in\fz$. Let $\fp$ be a parahoric subalgebra of 
$\fg_K:=K\ot\fg$ with pro-nilradical $\fp_n$ such that $\fp=\fz\op\fp_n$ as $\CC$-vector spaces.
By the last corollary in \cite{\KL, \S6}, there exists a non-empty subset $\fU$ of $\x+\fp_n$ (open in the 
power series topology) and $\s\in cl(W)$ such that for any $x\in\fU$, $x$ is regular
semisimple in a Cartan subalgebra of $\fg_K$ of type $\s$ (see \cite{\KL,\S1,\S6}). Note that $\s$ does not
depend on the choice of $\fU$. We expect that it does not depend on the choice of $\fp$ and that $g\m\s$ is 
a map $G@>>>cl(W)$ whose fibres are exactly the strata of $G$. By the identification 4.1(c) this induces an 
injective map $\un{cl}(W)@>>>cl(W)$ whose image is expected to be the subset $cl_0(W)$ in 4.10 and whose 
composition with the obvious map $cl(W)@>>>\un{cl}(W)$ is expected to be the identity map of $\un{cl}(W)$.

\enddocument